\documentstyle[12pt]{article}

\def\IA{{\cal I}_{\cal A}}
\def\IAs{{\cal I}_{{\cal A}_s}}

\def\pd{\partial}

\def\D{{\Delta}}

\def\Ds{{\Delta^*}}

\def\w{{\omega}}
\def\({ \left(  }
\def\){ \right) }
\def\cdt{{\hskip-0.02cm\cdot\hskip-0.02cm}}

\def\bn#1{\bar\nu_{#1}^*}

\def\ns#1{\nu_{#1}^*}
\def\l#1{l^{(#1)}}
\def\ta#1{\theta_{a_#1}}
\def\da#1{{\pd \; \over \pd {a_{#1}}}}

\def\en#1{{e_{\bn{#1}}}}

\def\X#1#2(#3)#4#5{ {$X_{#1}^#2(#3)_{#5}^{#4}$} }

\def\nl{n \cdt l}
\def\Jl{\hat J \cdt l}

 \renewcommand{\tilde}{\widetilde}

 \newcommand{\la}{\langle}
 \newcommand{\ra}{\rangle}
\font\sectionfont=cmbx12 at 12pt

\def\section#1{\vskip 1.5truepc\centerline{\hbox {{\sectionfont #1}}}
\vskip 1truepc\noindent\stepcounter{section}}
\font\authorfont=cmti12 at 14pt
\font\titlefont=cmbx12 at 16pt

\begin{document}
\setcounter{page}{1}

\def\sec#1{\vskip 1.5truepc\centerline{\hbox {{\sectionfont
#1}}}\vskip1truepc\noindent\noindent\stepcounter{section}}

\def\newsubsec#1#2{\vskip .3pc plus1pt minus 1pt\noindent {\bf #1}
{\rm #2}}

\newcommand{\eps}{\varepsilon}
\newcommand{\alp}{\alpha}
\newcommand{\gam}{\gamma}
\newcommand{\Aut}{{\rm Aut}\,}
\newcommand{\Out}{{\rm Out}\,}
\newcommand{\Mult}{{\rm Mult}\,}
\newcommand{\Inn}{{\rm Inn}\,}
\newcommand{\Irr}{{\rm Irr}\,}
\newcommand{\IBr}{{\rm IBr}}
\newcommand{\Ker}{{\rm Ker}\,}
\renewcommand{\Im}{{\rm Im}\,}
\newcommand{\Ind}{{\rm Ind}}
\newcommand{\diag}{{\rm diag}}
\newcommand{\soc}{{\rm soc}\,}
\newcommand{\End}{{\rm End}\,}
\newcommand{\sol}{{\rm sol}}
\newcommand{\Hom}{{\rm Hom}\,}
\newcommand{\rank}{{\rm rank}\,}
\newcommand{\Syl}{{\rm Syl}\,}
\newcommand{\Tr}{{\rm Tr}\,}
\newcommand{\tr}{{\rm tr}\,}
\newcommand{\Gal}{{\it Gal}\,}
\newcommand{\Spec}{{\rm Spec}\,}
\newcommand{\ad}{{\rm ad\,}}
\newcommand{\Sym}{{\rm Sym}\,}
\newcommand{\Char}{{\rm char}\,}
\newcommand{\pr}{{\rm pr}}
\newcommand{\Rad}{{\rm Rad}\,}
\newcommand{\abel}{{\rm abel}}
\newcommand{\codim}{{\rm codim}}
\newcommand{\CC}{{\Bbb C}}
\newcommand{\RR}{{\Bbb R}}
\newcommand{\QQ}{{\Bbb Q}}
\newcommand{\ZZ}{{\Bbb Z}}
\newcommand{\GG}{{\Bbb G}}
\newcommand{\SSS}{{\Bbb S}}
\newcommand{\AAA}{{\Bbb A}}
\newcommand{\FF}{{\Bbb F}}
\newcommand{\PP}{{\Bbb P}}
\newcommand{\HH}{{\Bbb H}}
\newcommand{\NN}{{\Bbb N}}
\newcommand{\UU}{{\Bbb U}}
\newcommand{\KK}{{\Bbb K}}
\newcommand{\LL}{{\Bbb L}}

 \catcode`\@=11
\font\twelvemsb=msbm10 scaled 1200
\font\tenmsb=msbm10
\font\ninemsb=msbm7 scaled 1200
\newfam\msbfam
\textfont\msbfam=\twelvemsb  \scriptfont\msbfam=\tenmsb
  \scriptscriptfont\msbfam=\ninemsb
\def\msb@{\hexnumber@\msbfam}
\def\Bbb{\relax\ifmmode\let\next\Bbb@\else
 \def\next{\errmessage{Use \string\Bbb\space only in math
mode}}\fi\next}
\def\Bbb@#1{{\Bbb@@{#1}}}
\def\Bbb@@#1{\fam\msbfam#1}
\catcode`\@=12

\renewcommand{\thefootnote}{\fnsymbol{footnote}}

\centerline{\titlefont  GKZ Systems, Gr\"obner Fans and Moduli Spaces}
\centerline{\titlefont   {of Calabi-Yau Hypersurfaces}}
\vskip 2pc
\centerline{\authorfont  S. Hosono }

\date{}

\vskip 3pc

\centerline{\sectionfont  Abstract}
\vskip 8pt
We present a detailed analysis of the GKZ(Gel'fand, Kapranov and Zelevinski) 
hypergeometric systems in the context of mirror symmetry of 
Calabi-Yau hypersurfaces in toric varieties. As an application we will derive 
a concise formula for the prepotential about large complex structure 
limits.

\section{1. Introduction} 
Mirror symmetry of Calabi-Yau manifolds 
has been playing a central role in revealing non-perturbative 
aspects of the type II string vacua, i.e., the moduli spaces for 
a family of Calabi-Yau manifolds. Since the success in determining the 
quantum geometry on the IIA moduli space made by Candelas et al\cite{CdGP} 
in 1991, there have been many progresses and a lot of 
communications between physics and mathematics on this topics\cite{GY}. 

In this article, we will be concerned with the mirror 
symmetry of Calabi-Yau hypersurfaces in toric varieties. 
In this case, the mirror symmetry may be traced to a rather 
combinatorial properties of the reflexive polytopes which determines 
the ambient toric varieties due to ref.\cite{Bat1}. 
Furthermore since the period integrals of Calabi-Yau hypersurfaces turn out 
to satisfy the hypergeometric differential equation, ${\cal A}$-hypergeometric 
system, introduced by Gel'fand, Kapranov and Zelevinski (GKZ), 
we can study in great detail the moduli spaces of Calabi-Yau hypersurfaces. 
Based on the analysis of GKZ-hypergeometric system in our context, we 
will derive a closed formula for the prepotential, which defines the special 
K\"ahler geometry on the moduli spaces. 

In section 2, we will review the mirror symmetry of Calabi-Yau 
hypersurfaces in toric varieties. This is meant to fix our notations as 
well as to introduce the mirror symmetry due to Batyrev. In section 3, 
we will introduce GKZ-hypergeometric system ($\Ds$-hypergeometric system) 
as an infinite set of differential equations satisfied by period integrals 
and summarize basic results following \cite{GKZ1}.  We also define the extended 
$\Ds$-hypergeometric system incorporating the automorphisms of the toric 
varieties. We will remark that the $\Ds$-hypergeometric system in our 
context is resonant in general. In section 4, we will review basic  
properties of the toric ideal and the Gr\"obner fan as an equivalence 
classes of the term orders in the toric ideal. We will use the Gr\"obner  
fan to compactify the space of the variables in the $\Ds$-hypergeometric 
system, and propose it as a natural compactification of the corresponding 
family of Calabi-Yau hypersurfaces. In section 5, we will prove general 
existence of the so-called large complex structure limits, at which 
the monodromy becomes maximally unipotent\cite{Mor}. We will also present 
a general formula for the local solutions about these points. In the 
final section, we will derive a closed formula for the prepotential, 
which is valid about a large complex structure limit for arbitrary 
Calabi-Yau hypersurfaces in toric varieties. Our formula determines 
the special K\"ahler geometry about a large complex structure limit 
as well as the quantum corrected Yukawa coupling. Claim 5.8, Claim 5.11, 
and Claim 6.8 in the last two sections are meant to state those results 
that are verified in explicit calculations by many examples without 
general proofs. 

All the results except Prop.6.7 for the prepotential in the 
final section have already reported in refs.\cite{HKTY1}\cite{HKTY2}
\cite{HLY1}\cite{HLY2}. 

\vspace{0.3cm}\noindent
{\bf Acknowledgments.} 
This article is based on the joint works 
with A. Klemm, B.H. Lian, S. Theisen and S.-T. Yau. 
The author would like to thank them for their enjoyable collaboration. 
He is also grateful to express his thanks to the organizers of the 
Taniguchi symposium, as well as the Taniguchi Foundation, 
where he had valuable discussions with many participants.  
He also express his thanks to L. Borisov and the referee for their 
valuable comments on the first version of this article.

\section{2. Mirror Symmetry of Calabi-Yau Hypersurfaces}

In this section, we will summarize mirror symmetry of 
Calabi-Yau hypersurfaces in toric varieties due to 
Batyrev.  We refer the paper\cite{Bat1} for details. 

Let $M\cong \ZZ^d$ be a lattice of rank $d$ and $N$ be its 
dual. We denote the dual pairing $M\times N \rightarrow \ZZ$ by 
$\la  , \ra $. A (convex) polytope $\D$ is a convex hull of a finite set 
of points in $M_{\RR}:=M\otimes\RR$. In the following, we assume $\D$ 
contains the origin in its interior. 
The polar dual $\Ds\subset N_{\RR}$ is defined by
\begin{equation}
\Ds=\{ x \in N_\RR \; \vert\; \la x,y\ra  \geq -1 ,\; y\in \D \;\}  \;\;.
\label{eqn: Ds}
\end{equation}

\vspace{0.3cm}\noindent
{\bf Definition 2.1.} 
A polytope $\D$ is called reflexive if it is a convex hull of 
a finite set of integral points in $M_\RR$ and contains only the 
origin in its interior. 

\vspace{0.3cm}\noindent
{\bf Proposition 2.2.} 
When a polytope $\D$ is reflexive, its dual $\Ds$ is also reflexive. 

\vspace{0.2cm}
Since $(\Ds)^*=\D$, reflexive polytopes come with a pair $(\D,\Ds)$. 
The following descriptions about $\Ds$ with $N$ equally apply to 
$\D$ with $M$ by symmetry. 

\vspace{0.3cm}\noindent
{\bf Definition 2.3.} 
A maximal triangulation $T_o$ of $\Ds$ is a simplicial decomposition 
of $\Ds$ with properties; 1) every $d$-simplex contains the origin 
as its vertex, 2) 0-simplices consist of all integral points of $\Ds$. 

\vspace{0.2cm}
For a maximal triangulation $T_o$ of $\Ds$, we consider a complete fan 
$\Sigma(\Ds,T_o)$ over the triangulation $T_o$ in $N_\RR$.  
Associated to the data $(\Sigma(\Ds,T_o),$ $N)$ 
we consider a toric variety $\PP_{\Sigma(\Ds,T_o)}$ \cite{Oda}\cite{Ful}. 
Due to the property that $\Ds$ is reflexive, we have 

\vspace{0.3cm}\noindent
{\bf Proposition 2.4.} (Prop.2.2.19 in \cite{Bat1}) {\it  
$\PP_{\Sigma(\Ds,T_o)}$ is a projective variety for at least 
one maximal triangulation with its anti-canonical class  
$-K=\displaystyle{ \sum_{\rho \in N\cap \Ds \setminus \{ 0\} } } D_\rho$ 
ample. }

\vspace{0.3cm} \noindent
{\bf Note.} In \cite{Bat1}, the maximal triangulations with the property 
in this proposition are called {\it projective}.  
In case of $d\leq 4$, we can observe widely that 
every maximal triangulation is projective. More generally we observe 
that every triangulation of a reflexive polytope 
is {\it regular} which generalize projective(, see right after 
eq.(\ref{eqn:convex}) for the definition). 
For a restricted class of reflexive polytopes 
(the type I or II in the following classification), 
it has been proved (Th.4.10 in \cite{HLY2}) that 
every nonsingular maximal triangulation is projective, 
see also Remark after Th.2.5. 
In the following, we will write a projective toric variety by 
$\PP_{\Sigma(\Ds,T_o)}$ choosing a projective maximal triangulation 
$T_o$ of $\Ds$. 

\vspace{0.2cm}
Let us fix a basis $\{ n_1,\cdots,n_d\}$ of $N$ and denote its 
dual basis by $\{m_1,\cdots,m_d\}$. 
With respect to this basis, we denote the coordinate ring of 
the torus $T_N:=\Hom_\ZZ(M,\CC^*) \subset \PP_{\Sigma(\Ds,T_o)}$ 
by $\CC[Y_1^{\pm1},\cdots,Y_d^{\pm1}]$ with the generators 
$Y_k={\rm{\bf e}}(m_k):T_N\rightarrow \CC^*$ defined by 
{\bf e}$(m_k)(t)=t(m_k)$. Consider a Laurent 
polynomial $f_\D = \sum_{\nu\in \D\cap M} c_{\nu} Y^\nu$ with 
complex coefficients $c_\nu$. We denote by $X_\D$ the Zariski closure of 
the zero locus $(f_\D=0)$ in $\PP_{\Sigma(\Ds,T_o)}$ for generic $c_\nu$'s. 
Similarly, we consider a projective toric variety $\PP_{\Sigma(\D,T_o)}$ 
associated to a projective maximal triangulation $T_o$ of $\D$, and 
denote the coordinate ring of $T_M:=\Hom_\ZZ(N,\CC^*) 
\subset \PP_{\Sigma(\D,T_o)}$ by $\CC[X_1^{\pm1},\cdots,X_d^{\pm1}]$ with 
$X_k={\rm{\bf e}}(n_k)$.

\vspace{0.3cm}\noindent
{\bf Theorem 2.5.} (Th.4.2.2, Corollary 4.2.3, Th.4.4.3 in \cite{Bat1}) 
{\it Let $(\D,\Ds)$ be a pair of reflexive polytopes in 
dimensions $d\leq4$ (, in $M_\RR$ and $N_\RR$, respectively). Then; 
1) generic hypersurfaces $X_\D\subset \PP_{\Sigma(\Ds,T_o)}$ and 
$X_\Ds \subset \PP_{\Sigma(\D,T_o)}$ 
define smooth Calabi-Yau manifolds, 
2) these two hypersurfaces are mirror symmetric in their 
Hodge numbers, i.e., $h^{1,1}(X_\D)=h^{d-2,1}(X_\Ds) , 
\;$ $ h^{d-2,1}(X_\D)= h^{1,1}(X_\Ds)$. }

\vspace{0.3cm}
\noindent
{\bf Remark.} 
Depending on the toric data of the reflexive polytopes, 
the ambient spaces have (Gorenstein) singularities (Prop.2.2.2 in \cite{Bat1}) 
in general. We call a maximal triangulation is {\it nonsingular} 
if its simplices of maximal dimensions consists of unit simplices, i.e.,
simplices with unit volume. It is easy to deduce that the toric variety 
is nonsingular if the maximal triangulation is so. 
Now we classify the reflexive polytopes into the following 
three types:  
\par\noindent
$\bullet$ type I;
the polytope 
has no integral point in the interior of all codimension one faces, 
and has a nonsingular maximal triangulation,
\par\noindent 
$\bullet$ type II;
the polytope has at least one integral point in the interior of some 
codimension one face, and has a nonsingular maximal triangulation,
\par\noindent
$\bullet$ type III;
the polytope does not have a nonsingular maximal triangulation. 

In the following, we always consider a nonsingular maximal triangulation 
$T_o$ for the polytopes of type I and II. Then the toric varieties 
$\PP_{\Sigma(\Ds,T_o)}$ are projective and nonsingular for both 
the polytopes $\Ds$ of type I and II (Th.4.10 in \cite{HLY2}), 
however we distinguish these two 
because of the difference in the {\it root system} for their 
dual polytopes $\D$; 
\begin{equation}
\begin{array}{rcl}
R(\D,M)=\{ \alpha \in \Ds\cap N \;\vert\; 
 && \makebox[-1em]{}
  \exists m_\alpha \in \D\cap M \; s.t. 
  \la  m_\alpha, \alpha \ra =-1 \;   \\ 
 &&\makebox[-1em]{} {\rm and}\;  
  \la  m , \alpha \ra  \geq 0 \;(m\not= m_\alpha \in \D\cap M) \;\; \}  \\
\end{array}
\label{eqn: rootD}
\end{equation}
The root system determines 
the automorphisms of the toric variety $\PP_{\Sigma(\D,T_o)}$ infinitesimally 
due to the following result, which we will utilize in the next section;

\vspace{0.3cm}\noindent
{\bf Proposition 2.6.} (Prop.3.13 in \cite{Oda}) {\it 
For a nonsingular toric variety $\PP_{\Sigma(\D,T_o)}$, we have a direct sum 
decomposition via the root system $R(\D,M)$;  
\begin{equation}
\begin{array}{rcl}
{\rm Lie}(\Aut(\PP_{\Sigma(\D,T_o)})) &=& 
H^0(\PP_{\Sigma(\D,T_o)},\Theta_{\PP_{\Sigma(\D,T_o)}}) \\
&=& {\rm Lie}(T_M) \oplus \( 
\oplus_{\alpha\in R(\D,M)}
\CC {\rm{\bf e}}(\alpha)\delta_{m_{\alpha}}\) \\
\end{array}
\label{eqn: lie}
\end{equation}
where $\delta_m \; (m \in M)$ is the derivation on $T_M$ defined by 
$\delta_m {\rm{\bf e}}(n):=\la  m,n \ra {\rm{\bf e}}(n)$. }

\newpage
\section{3. Resonance in GKZ Hypergeometric System}

We consider a family of Calabi-Yau hypersurfaces $X_{\Ds}(a) 
\subset \PP_{\Sigma(\D,T_o)}$ varying the coefficients $a_{\ns{{}}}$ in the 
defining equation $f_\Ds(a)=\sum_{\ns{{}}\in \Ds\cap N} a_{\ns{}} X^{\ns{}}$. 
By this polynomial deformation, we describe the complex structure 
deformation of $X_\Ds$. This deformation space is mapped to that 
of (complexified) K\"ahler class of $X_\D\subset \PP_{\Sigma(\Ds,T_o)}$ 
under the mirror symmetry. According to the local Torelli theorem \cite{BG}, 
we can introduce a local coordinate on the moduli space in terms of 
period integrals. In case of hypersurfaces in toric varieties, we 
have one canonical period integral \cite{Bat2}\cite{BC}
\begin{equation}
\Pi(a)={1\over (2\pi i)^d} \int_{C_0} {1\over f_\Ds(X,a)} 
\prod_{i=1}^{d} {d X_i \over X_i} \;\;,
\label{eqn: Pi}
\end{equation}
with the cycle $C_0=\{ \vert X_1\vert=\cdots=\vert X_d\vert=1\}$ in $T_M$. 
Here we study the differential equation satisfied by (\ref{eqn: Pi}).

\newsubsec{(3-1) Extended GKZ hypergeometric system } 
Let ${\cal A}=\{ \bar\chi_0,\cdots,\bar\chi_p\}$ be a finite 
set of integral points in $\{1\}\times\RR^n\subset\RR^{n+1}$. We assume
the vectors $\bar\chi_0,\cdots,\bar\chi_p$ span $\RR^{n+1}$. 

\vspace{0.3cm}\noindent
{\bf Definition 3.1.} 
Consider the lattice of relations among the set ${\cal A}$, 
\begin{equation}
L=\{ \; (l_0,\cdots,l_p)\in \ZZ^{p+1}\; \vert\; 
\sum_{i=0}^p l_i\bar\chi_{i,j}=0 \;, (j=1,\cdots,n+1) \;\} \;,
\label{eqn: LA}
\end{equation}
where $\bar\chi_{i,j}$ represents the $j$-th component of 
the vector $\bar\chi_i$.  {\it ${\cal A}$-hypergeometric 
system with exponent $\beta \in \CC^{n+1}$ } is a system of 
differential equations for a complex function $\Psi(a)$ 
on $\CC^{\cal A}$;
\begin{eqnarray}
{\cal D}_l \Psi(a)&=&\{
\prod_{l_i>0}\({\pd \; \over \pd a_i}\)^{l_i} -
\prod_{l_i<0}\({\pd \; \over \pd a_i}\)^{-l_i} \} \Psi(a)=0 
\;\; (l \in L ) \\
{\cal Z}\Psi(a)&=&\{ \sum_{i=0}^p \bar\chi_i a_i {\pd \; \over \pd a_i} 
-\beta \} \Psi(a) =0  \;\; .
\label{eqn: gkz}
\end{eqnarray}

\vspace{0.3cm}\noindent
{\bf Proposition 3.2.} (\cite{Bat2}) {\it 
The period integral (\ref{eqn: Pi}) satisfies 
the ${\cal A}$-hyper-\break geometric system 
with ${\cal A}=\{1\}\times (\Ds\cap N)$ and 
$\beta=(-1)\times\vec 0$. We call this hypergeometric system as 
$\Ds$-hypergeometric system. }

\vspace{0.2cm}
By direct evaluation of the action of ${\cal D}_l$ and ${\cal Z}$ 
on the period integral (\ref{eqn: Pi}), we obtain this proposition. 
Here we consider the meaning of the linear operator ${\cal Z}$. 
The first component of this vector operator is exactly the 
Euler operator, and just says that the period integral has homogeneous 
degree $-1$ as a function of $a_i$'s. For the other components, it is 
easy to deduce that these come from the invariance of the period 
integral under the torus actions, which act infinitesimally on 
the coordinate $X_k={\rm{\bf e}}(n_k)$ by 
$\delta_m X_k = \la  m, n_k \ra  X_k$. It is now clear that these actions 
should be considered for all elements in 
${\rm Lie}(\Aut (\PP_{\Sigma(\D,T_o)}))$. 
Then we may write the invariance of the period integral under the 
infinitesimal action of $\xi \in {\rm Lie}(\Aut (\PP_{\Sigma(\D,T_o)}))$ 
by a linear differential operator ${\cal Z}_\xi$ acting on $\Pi(a)$ through  
\begin{equation}
{\cal Z}_{\xi} \Pi(a)=
\int_{C_0} \xi  
\( {1\over f_\Ds(X,a)} \) 
\prod_{i=1}^{d} {d X_i \over X_i}  = 0 \;\;.
\label{eqn: Zxi}
\end{equation}
For explicit forms of the operators ${\cal Z}_\xi$, we refer to the examples 
given in p.541 of \cite{HLY1}.

\vspace{0.3cm}\noindent
{\bf Proposition 3.3.} ((2.13) in \cite{HLY1}) 
{\it The period integral $\Pi(a)$ satisfies 
\begin{equation}
\begin{array}{crl}
&& {\cal D}_l \Pi(a) =0 \; (l\in L) \;,\; \\
&& \\
&& {\cal Z}_E \Pi(a)=0 \;,\;
{\cal Z}_\xi \Pi(a)= 0 \;\; 
(\xi \in {\rm Lie}(\Aut (\PP_{\Sigma(\D,T_o)}))), \\
\end{array}
\label{eqn: extGKZ}
\end{equation}
where we denote the Euler operator by 
${\cal Z}_E=\sum_{i=0}^p a_i{\pd\; \over \pd a_i} +1$. }

\vspace{0.2cm}
We call this system as {\it extended GKZ-hypergeometric system} or 
{\it extended $\Ds$-hypergeometric system}.  By Prop.2.6, it is 
clear that this extended system reduces to the GKZ system if the 
polytope $\Ds$ is of type I. In the following, we take an approach to 
study mainly the $\Ds$-hypergeometric system because the set of the solutions 
of the extended $\Ds$-hypergeometric system can be found in that of the 
$\Ds$-hypergeometric system.

\newsubsec{(3-2) Convergent series solutions } 
Here we summarize general results in \cite{GKZ1} about the 
convergent series 
solution of the ${\cal A}$-hypergeometric system with exponent $\beta$. 
This is to introduce the notion of the secondary fan as well as to 
fix our conventions and notations. Since our interest is in the period 
integrals, we assume ${\cal A}=\{1\}\times (\Ds\cap N)$ and
$\beta=(-1)\times \vec 0$. Hereafter we write the integral points
explicitly by $\Ds\cap N=\{\ns{0},\cdots,\ns{p}\} \;(\ns{0}\equiv\vec
0)$ and $\bn{i}:=1\times\ns{i} \; (i=0,\cdots,p)$. 
 
We start with a formal solution of the ${\cal A}$-hypergeometric system
with exponent $\beta$ given by
\begin{equation}
\Pi(a,\gamma)=\sum_{l\in L} {1 \over 
\prod_{0\leq i\leq p} \Gamma(l_i+\gamma_i+1)}a^{l+\gamma} 
\;\;,
\label{eqn: fsol}
\end{equation}
where $\beta=\sum_i \gamma_i \bn{i}$. Now define an affine space 
$\Phi(\beta):=\{ \gamma\in \RR^{p+1}\vert \beta=\sum \gamma_i\bn{i} \}$. 
A subset $I\subset \{ 0,\cdots,p\}$ is called a base if $\bn{I}:=\{\bn{i}\vert 
i\in I\}$ form a basis of $\RR^{d+1}$. Given a base $I$, we may solve 
$\sum_{j\in I}\gamma_j \bn{j} = \beta -\sum_{j\not\in I}\gamma_j\bn{j}$ 
for $\gamma_j \;(j\in I)$ and define $\Phi_{\ZZ}(\beta,I)=\{ 
\gamma\in \Phi(\beta) \;\vert\; \gamma_j \in \ZZ \;(j\not\in I) \}$. 
We choose an integral basis $A=\{ \l1, \cdots, \l{p-d} \}$ of the lattice
$L$, and define $\Phi^A_{\ZZ}(\beta,I)=\{ \gamma\in \Phi_{\ZZ}(\beta,I)
\;\vert\; \gamma_j=\sum_{k=1}^{p-d} \lambda_k \l{k}_j \; (0\leq \lambda_k<1,
\; j\not\in I) \}$. Then $\Phi^A_{\ZZ}(\beta,I)$ provides 
a set of representatives of the quotient $\Phi_\ZZ(\beta,I)/L$ 
and kills the invariance $\gamma \rightarrow \gamma+v \; (v\in L)$ in the 
formal solution (\ref{eqn: fsol}). 

\vspace{0.3cm}\noindent
{\bf Definition 3.4.} For a base $I$, define a cone in $L_\RR=L\otimes\RR$
by ${\cal K}({\cal A},I)=\{ l \in L_\RR \;\vert\; l_i\geq0\;(i\not\in I)\}$. 
A $\ZZ$-basis $A\subset L$ is said {\it compatible} 
with a base $I$ if it generates a cone that contains ${\cal K}({\cal
A},I)$. 

\vspace{0.3cm}\noindent
{\bf Theorem 3.5.} (Prop.1 in \cite{GKZ1}) {\it 
Fix a base $I$ and choose a $\ZZ$-basis $A=\{\l1,\cdots,\l{p-d}\}$ 
compatible with  
it. Then the formal solution (\ref{eqn: fsol}) 
takes the form $\Pi(a,\gamma)=a^\gamma
\sum_{m\in \ZZ^{p-d}_{\geq0}}c_m(\gamma) x^m$ for each $\gamma\in 
\Phi_\ZZ^A(\beta,I)$ with $x_k=a^{\l{k}}$. This powerseries converges 
for sufficiently small $\vert x_k \vert$. }

\vspace{0.3cm}\noindent
{\bf Remark.} 
The coefficient $c_m(\gamma)$ is given explicitly by 
$$
c_m(\gamma)={1\over \prod_{i=0}^p \Gamma(\sum_k m_k\l{k}_i + \gamma_i+1)} \;.
$$
For some index $i$ of the base $I$ in the above theorem it can happen that 
$\sum m_k\l{k}_i+\gamma_i+1$ is non-positive for all $m \in
\ZZ_{\geq0}^{p-d}$ , which means we have the trivial solution 
$\Pi(a,\gamma)\equiv0$. In this case, we multiply an infinite number 
$\Gamma(\gamma_i+1)$ to obtain 
nonzero powerseries, i.e., 
$
{\Gamma(\gamma_i+1) \over \Gamma(\sum m_k\l{k}_i+\gamma_i+1)} := 
\lim_{\varepsilon\rightarrow0}
{\Gamma(\gamma_i+1+\varepsilon) \over 
\Gamma(\sum m_k\l{k}_i+\gamma_i+1+\varepsilon)} .
$
 
\vspace{0.3cm}\noindent
{\bf Definition 3.6.} 
Consider $P={\rm Conv.}\( \{ 0, \bn{0},\cdots,\bn{p} \}\)$ in
$\RR^{d+1}$. A collection of bases $T=\{I\}$ is called a triangulation
of $P$ if $\cup_{I\in T}\la \bn{I}\ra  =P$ for simplices
$\la \bn{I}\ra ={\rm Conv.}\(\{0\}\cup\bn{I}\)$ $(I\in T)$, 
and $\la \bn{I_1}\ra \cap\la \bn{I_2}\ra \;(I_1,I_2\in T)$ is 
a lower dimensional common face.  

\vspace{0.3cm}\noindent
{\bf Note.}
Since $(d+1)$-simplices in $P$ are in one-to-one correspondence to 
$d$-simplices in $\Ds$, we identify a triangulation of $T$ of $P$ 
with its corresponding triangulation of $\Ds$. We call 
a triangulation $T$ of $P$ is maximal if it corresponds to a maximal 
triangulation of $\Ds$(Def.2.3). 

\vspace{0.2cm}
For a base $I$ and a point $\eta\in \RR^{p+1}$, we consider a linear
function $h_{I,\eta}$ on $\RR^{p+1}$ by $h_{I,\eta}(\bn{i})=\eta_i \;(i\in I)$.
We define ${\cal C}({\cal A},I)=\{ \eta\in \RR^{p+1} \;\vert\;
h_{I,\eta}\;(\bn{i}) \leq \eta_i \;(i\not\in I)\}$ and ${\cal C}({\cal
A},T):=\cap_{I\in T}{\cal C}({\cal A},I)$ for a triangulation $T$.
Then it is easy to see that ${\cal C}({\cal A},T)$ consists of 
$\eta \in \RR^{p+1}$ for which we have a convex piecewise 
linear function $h_{T,\eta}$ on $T$ determined by 
$h_{T,\eta}\vert_{\la \bn{I}\ra}= h_{I,\eta} \; (I\in T)$ and satisfies 
$h_{T,\eta}(\bn{i})\leq \eta_i$ for $\bn{i}$ 
not a vertex of $T$, i.e., 
\begin{equation}
\begin{array}{crl}
{\cal C}({\cal A},T)=\{ \eta \in \RR^{p+1}\; \vert \;\; 
&& \makebox[-2em]{}
   h_{I_1,\eta}(v)\leq h_{I_2,\eta}(v) \;\; 
   (v\in \la \bn{I_2}\ra , \; I_1,I_2\in T), \\
&& \makebox[-2em]{}
   h_{T,\eta}(\bn{i})\leq \eta_i \;\;
   (\bn{i} \;\,{\rm is} \;\, {\rm not}\;\, {\rm a} \;\, 
{\rm vertex}\;\,{\rm of}\; T)\;  \}\\
\end{array}
\label{eqn:convex}
\end{equation}
A triangulation is called {\it regular} if ${\cal C}({\cal A},T)$ has 
interior points. 
We say a $\ZZ$-basis $A\subset L$ is compatible with a triangulation 
$T$ if it is compatible with all bases $I$ in $T$. 

\vspace{0.3cm}\noindent
{\bf Proposition 3.7.} (Prop.5 in \cite{GKZ1}) {\it 
For every regular triangulation $T$, there exists (infinitely many) 
$\ZZ$-basis of $L$ compatible with $T$. }

\vspace{0.2cm}
The exponent $\beta$ is called {\it $T$-nonresonant} if the set 
$\Phi^A_\ZZ(\beta,I)\; (I\in T)$ are pairwise disjoint. 
We set $\Phi^A_\ZZ(\beta,T):=\cup_{I\in T}\Phi^A_\ZZ(\beta,I)$. 
We normalize the volume of the standard $(d+1)$-simplex to $1$. 

\vspace{0.3cm}\noindent
{\bf Theorem 3.8.} (Th.3 in \cite{GKZ1}) {\it 
For a regular triangulation $T$ of the polytope $P$, and a $\ZZ$-basis 
$A=\{\l{1},\cdots,\l{p-d}\}$ of $L$ compatible with $T$, there are integral 
powerseries in the variables $x_k=a^{\l{k}}$ for 
$a^{-\gamma}\Pi(a,\gamma) \;(\gamma\in \Phi^A_\ZZ(\beta,T))$, which 
converge for sufficiently small $\vert x_k\vert$. If the exponents 
$\beta$ is $T$-nonresonant, these series constitute $vol(P)$ linearly 
independent solutions. }

\vspace{0.3cm}\noindent
{\bf Remark.} 
In our case of $\Ds$-hypergeometric system with $\beta=(-1)\times\vec0$, 
we have one special element $\gamma=(-1,0,\cdots,0)$ in the set 
$\Phi_\ZZ(\beta,I)$ for any base $I$. If the polytope $\Ds$ is of 
type I or II in our classification, a nonsingular maximal 
triangulation $T$ of $P$ consists of those bases $I$ for which 
$\vert \det (\bar\nu_{j,i}^*)_{1\leq i\leq d+1,\; j\in I}\vert =1$. 
Because of this unimodularity, we have 
\begin{equation}
\Phi_\ZZ(\beta,I)=(-1,0,\cdots,0)+L \;\; , 
\label{eqn: phiL}
\end{equation}
and $\Phi_\ZZ^A(\beta,I)=\{ (-1,0,\cdots,0)\}$ for every base 
$I$ of the maximal triangulation and any $\ZZ$-basis $A$ compatible with it. 
Thus we encounter a ``maximally $T$-resonant'' situation. 

\vspace{0.3cm}\noindent
{\bf Definition 3.9.} 
For a regular triangulation $T$ and a $\ZZ$-basis 
$A=\{ \l{1},\cdots,$ $\l{p-d}\}$ 
compatible with $T$, we define a 
power series $w_0(x,\rho; A)=a_0\Pi(a,\gamma)$ 
(with $\gamma=\sum_{k=1}^{p-d}\rho_k\l{k}+(-1,0,\cdots,0)$) by
\begin{equation}
w_0(x,\rho;A)=\sum_{m\in \ZZ_{\geq0}^{p-d} }
{ \Gamma(-\sum_k (m_k+\rho_k)\l{k}_0+1) \over
\prod_{1\leq i\leq p}
\Gamma(\sum_k(m_k+\rho_k)\l{k}_i+1) } x^{m+\rho} \;\;,
\label{eqn: wnot}
\end{equation}
where  
$x_k=(-1)^{\l{k}_0}a^{\l{k}}$.

\vspace{0.3cm}\noindent
{\bf Remark.}  
Here we have applied our recepie of multiplying the constant 
$\Gamma(\gamma_0+1)$ to the formal solution $\Pi(a,\gamma)$. 
We adopt this definition because for a maximall triangulation $T_o$, 
we encounter the situation $\Pi(a,\gamma)\equiv 0$, namely, 
$\sum m_k \l{k}_0 +\gamma_0+1 \in \ZZ_{\leq0} \; 
(m\in \ZZ^{p-n}_{\geq0})$ for a $\ZZ$-basis 
$A=\{\l{1},\cdots,\l{p-n}\}$ compatible with $T_o$. 
(See Prop.4.8, Prop.4.9 and Th.4.10 in \cite{HLY2}).  
In general a basis $A$ compatible with a regular triangulation $T$ 
contains both the bases vectors $\l{k}$ with positive 0-th component 
and nonpositive 0-th component. Taking a value 
$\rho_0 \in \ZZ^{p-n}_{\geq0}$ under this situation cause infinity 
for some $m$ in the numerator of the coefficients of (\ref{eqn: wnot}). 
In this case we understand in our definition (\ref{eqn: wnot}) 
to take a limit $\rho \rightarrow \rho_0$ in a generic way. 
(If we encounter infinities under this limit in some coefficient, 
we go back to the original definition (\ref{eqn: fsol}). 
For a maximall triangulation $T_o$, we observe that this limit exists 
for all coefficients in (\ref{eqn: wnot}), see Claim5.8.)

\newsubsec{(3-3) Secondary fan}  
It is known that the set ${\cal C}({\cal A},T)$ is a closed 
polyhedral cone and that these cones cover $\RR^{p+1}$ when we vary the 
triangulations. Thus the set of these cones and their lower dimensional faces 
all together define a complete, polyhedral fan ${\cal F}({\cal A})$ 
called the {\it secondary fan}\cite{BFS}\cite{OP}.  
Let $\overline M = \ZZ \times M$ and 
$\overline 
N = \ZZ \times N$. We extend naturally the pairing $\la  ,\ra $ to that 
of $\overline M$ and $\overline N$. Consider the lattice 
$\overline {\cal M}:= \oplus_{\ns{} \in \Ds\cap N}
\ZZ \en{}\;(=\oplus_{i=0}^p \ZZ\en{i})$ and its dual 
$\overline {\cal N}=\Hom_\ZZ(\overline {\cal M},\ZZ)$. Then we have the 
following exact sequences
\begin{equation}
\begin{array}{crl}
&& \makebox[2.5cm]{} 
0 
\smash{\mathop{\longrightarrow}\limits^{}} \; \overline M \;
\smash{\mathop{\longrightarrow}\limits^{\quad {\bf A}\quad}} 
    \; \overline{\cal M} \;
\smash{\mathop{\longrightarrow}\limits^{\quad {\bf B}\quad }} 
    \;  \Xi(\overline M) \;
\smash{\mathop{\longrightarrow}\limits^{}}  0 \;\;,  \\ 
&& 
0  \smash{\mathop{\longleftarrow}\limits^{}}  \; {\rm Coker}{\bf A}^* \; 
\smash{\mathop{\longleftarrow}\limits^{}} \; \overline N \;
\smash{\mathop{\longleftarrow}\limits^{\quad {\bf A}^* \quad}} 
     \; \overline{\cal N} \;
\smash{\mathop{\longleftarrow}\limits^{\quad {\bf B}^* \quad}} 
     \; \Xi(\overline N) \;
\smash{\mathop{\longleftarrow}\limits^{}}  0 \;\;,  \\ 
\end{array}
\label{eqn: exacts}
\end{equation} 
where ${\bf A}(\bar m)=\sum_{i=0}^p\la \bar m,\bn{i}\ra \en{i} \;\; 
(\bar m\in \overline M)$ and ${\bf B}$ is the quotient. The dual is given by 
${\bf A}^*(\mu)=\sum_{i=0}^p \mu(\en{i})\bn{i} \;\; (\mu\in \overline N)$.
The pair $\{ {\cal B}, \Xi(\overline M) \}$ with 
${\cal B}:=\{ {\bf B}(\en{0}), \cdots, {\bf B}(\en{p}) \}$ is 
called {\it Gale transform} of a pair 
$\{ {\cal A}, $ $\overline N \}$. Under this general 
setting, let us consider a polyhedral cone in $\Xi(\overline M)_\RR$ 
\begin{equation}
{\cal C}'({\cal A},T)=\cap_{I\in T} \( \sum_{i\not\in I} 
\RR_{\geq0}{\bf B}(\en{i}) \)\;\;.
\end{equation} 

\vspace{0.3cm}\noindent
{\bf Proposition 3.10.}(Lemma 4.2 in \cite{BFS}) {\it 
The map ${\bf B}$ induces the following decomposition
\begin{equation}
{\cal C}({\cal A},T)={\rm Ker}({\bf B})\oplus {\cal C}'({\cal A},T) \;\;.
\label{eqn:decomp}
\end{equation} }

\vspace{0.3cm}
By definition, the cone ${\cal C}'({\cal A},T)$ is strongly convex. 
Using the above decomposition, we redefine the secondary fan to be 
\begin{equation}
{\cal F}({\cal A})=\{ {\cal C}'({\cal A},T) \vert \; T: 
{\rm regular} \;\; {\rm triangulation} \} .
\end{equation}
Now the secondary fan consists of strongly convex, polyhedral cones. 
If the polytope $\Ds$ is of type I or II, then the quotient $\Xi(\overline M)$ 
is torsion free and thus 
$({\cal F}({\cal A}),\Xi(\overline M))$ defines a toric 
variety. Even in the case of type III, we may consider the corresponding 
toric variety by simply projecting out the torsion part of $\Xi(\overline M)$. 
We will use this toric variety for the compactification of the moduli 
space in the next section. 

Now let us note that $\Xi(\overline N)\cong {\rm Ker}({\bf A}^*)$ is 
identified with the lattice $L$, and thus ${\cal K}({\cal A},T) \subset 
\Xi(\overline N)_\RR$. By definition of ${\cal C}'({\cal A},T)$, we may 
deduce 
\begin{equation}
{\cal K}({\cal A},T)={\cal C}'({\cal A},T)^\vee \;\;.
\label{eqn: Cdual}
\end{equation}
Since for a regular triangulation $T$, ${\cal C}'({\cal A},T)$ is a strongly 
convex polyhedral cone with interior points, the dual cone 
${\cal K}({\cal A},T)$ is also strongly convex polyhedral cone. 
Therefore we see that there are infinitely many $\ZZ$-basis of the 
lattice $L$ compatible with $T$ (Prop.3.7).

\section{4. Toric Ideal and Gr\"obner Fan}

In this section we will reduce the infinite set of 
operators ${\cal D}_l\; (l\in L)$ in our $\Ds$-hypergeometric system 
to a finite set. This will be related to the compactification 
problem of the moduli spaces.

\newsubsec{(4-1) Toric ideal and Gr\"obner fan} 
We write the operators ${\cal D}_l$ in (6) simply by 
${\cal D}_l=({\pd\;\over\pd a})^{l_+} - ({\pd\;\over\pd a})^{l_-}$ with 
$l=l_+-l_-$.  Keeping this form in mind we define 
{\it toric ideal} in a polynomial ring:

\vspace{0.3cm}\noindent
{\bf Definition 4.1.} Consider a polynomial ring 
$\CC[y]:=\CC[y_0,\cdots,y_p]$. 
Toric ideal $\IA$ is defined to be an ideal generated by {\it binomials} 
$y^{l_+}-y^{l_-} \; (l\in L)$,
\begin{equation}
\IA=\la  \; y^{l_+}-y^{l_-} \;\vert\; l\in L \; \ra  \;. 
\end{equation}

\vspace{0.3cm}
Toric ideal has been extensively studied in ref.\cite{Stu1}\cite{GKZ2}. 
Here we summarize relevant results for our purpose. 

A {\it term order} (monomial ordering) on $\CC[y]$ is a total order $\prec$ 
on the set of monomials $\{ y^\alpha \;\vert\; \alpha\in \ZZ_{\geq0}^{p+1}\}$ 
satisfying, 1) $y^\alpha \prec y^\beta$ implies $y^{\alpha+\gamma} \prec 
y^{\beta+\gamma}$ and 2) $1$ is the unique minimal element. 
When we have an ideal ${\cal I} \subset \CC[y]$ and fix a term order, 
we can speak of the leading term $LT_\prec(f)$ for 
every non-zero polynomial in ${\cal I}$. 
Then we define {\it initial ideal} of ${\cal I}$ by
\begin{equation}
\la LT_\prec ({\cal I})\ra =\la  LT_\prec (f) \;\vert\; f\in {\cal I}, f\not=0 \ra \;.
\end{equation}
A finite set ${\cal G}\subset {\cal I}$ is called {\it Gr\"obner basis} with 
respect to a term order $\prec$ if it generates the initial ideal; 
\begin{equation}
\la  LT_\prec ({\cal I})\ra =\la  LT_\prec(g) \;\vert\; g \in {\cal G} \ra  \;.
\label{eqn: groebner}
\end{equation}

\vspace{0.3cm}\noindent
{\bf Theorem 4.2.} (Th.1.2 in \cite{Stu2}) {\it 
For every ideal ${\cal I}\subset \CC[y]$, there are only finitely many 
distinct initial ideals.  }


\vspace{0.3cm}
We consider representing the term orders by weight vectors 
$\w = (w_0,\cdots,$ $\w_p)\in \RR^{p+1}$. For a polynomial $f=\sum_\alpha 
c_\alpha y^\alpha$, we define its {\it leading terms} $LT_\w(f)$ to 
be a sum of terms $c_\alpha y^\alpha$ whose weight $t_\w(y^\alpha)
:=\w_0\alpha_0+\cdots+\w_p\alpha_p$ is maximal. It is obvious that if 
the components of $\w \in \RR_{\geq0}^{p+1}$ are rationally independent, 
the weight determines a term order on $\CC[y]$.   When we fix an ideal 
${\cal I}\subset \CC[y]$, we may relax the condition for the weight $\w$ 
to be a term order; 
we say a weight $\w\in \RR^{p+1}$ defines a term order 
of ${\cal I}$ if $\la LT_\w({\cal I})\ra =\la LT_\prec({\cal I})\ra $ for 
some term order $\prec$. The following proposition provides a 'converse' 
statement,

\vspace{0.3cm}\noindent
{\bf Proposition 4.3.} (Prop.1.11 in \cite{Stu2}) {\it 
For any term order $\prec$, there exists a weight 
$\w\in \RR_{\geq0}^{p+1}$ such that $\la LT_w({\cal I})\ra 
=\la LT_\prec({\cal I})\ra $. }

\vspace{0.3cm}
Now {\it Gr\"obner region} is defined to be a set 
\begin{equation}
GR({\cal I})=\{ \w\in \RR^{p+1} \vert 
\la  LT_\w({\cal I})\ra =\la  LT_{\w'}({\cal I})\ra  \;
{\rm for} \;{\rm some} \; \w'\in \RR_{\geq0}^{p+1} \} \; .
\end{equation}

\vspace{0.3cm}\noindent
{\bf Proposition 4.4.} (Prop.1.12 in \cite{Stu2}) {\it 
If an ideal ${\cal I}\in \CC[y]$ is a homogeneous ideal with some grading 
$deg(y_i)=d_i >0$, then $GR({\cal I})=\RR^{p+1}$. }

\vspace{0.3cm}
Since the toric ideal $\IA$ is homogeneous ideal with 
$deg(y_0)=\cdots=$ $deg(y_p)$ $=1$, we see $GR(\IA)=\RR^{p+1}$. 
For a term order $\w$ of $\IA$, we define
\begin{equation}
{\cal C}(\IA,\w)=\{ \w' \in \RR^{p+1} \;\vert\; 
\la  LT_\w(\IA)\ra = \la LT_{\w'}(\IA)\ra  \;\} \;.
\label{eq:Cw}
\end{equation}
It is known that this set constitutes an open, convex, polyhedral cone 
in $\RR^{p+1}$ (Prop.2.1 in \cite{Stu1}). In the following, we mean by 
${\cal C}(\IA, \w)$ the closure of the set (\ref{eq:Cw}). 
Then due to Th.4.2 and Prop.4.4, the collection $\{ {\cal C}(\IA, \w) \}$ is
finite and defines a complete polyhedral fan ${\cal F}(\IA)$ in $\RR^{p+1}$, 
called the {\it Gr\"obner fan}.

\newsubsec{(4-2) Indicial ideal and compactification of $\Hom_\ZZ(L,\CC^*)$ } 
In the previous section, we called a triangulation $T$ of the
polytope $P$ regular if the cone ${\cal C}({\cal A},T)$ has interior
points. Here we characterize the regular triangulation in a geometrical
way. To this aim let us first consider a polytope $P_\w:={\rm Conv.}\(\{ 
\w_0\times\ns{0},\cdots, \w_p\times\ns{p}\} \)$ in $\RR^{d+1}$ for 
a weight $\w \in \RR^{p+1}$. If we
project a polytope $P_\w$ to $1\times \RR^d$, then we have the polytope 
$1\times \Ds$. Thus we may regard the weight $\w$ giving a hight to each
vertex of $1\times\Ds$. For generic weight $\w$, the 'lower' faces of
the polytope $P_\w$ consist of simplices and define, under the
projection, a simplicial decomposition of $\Ds$ and thus induce a 
triangulation $T_\w$. The regular triangulation of the polytope $P$ is a
triangulation $T_\w$ obtained for some weight $\w$ in this way 
(see Def.5.3 of \cite{Zie} for more details). 
It is not difficult to see the relation of the polytope $P_\w$ to  
the piecewise linear function $h_{T,\eta}$ in (\ref{eqn:convex}) 
with $\eta=\w$. 

Given a (regular) triangulation $T$ of the polytope $P$, the
{\it Stanley-Reisner ideal} $SR_T$ in $\CC[y]$ is defined to be the ideal
generated by all monomials $y_{i_1}\cdots y_{i_k}$ for which the
vertices $\bn{i_1},\cdots,\bn{i_k}$ do not make a simplex in $T$. 
The following theorem is due to Sturmfels:

\vspace{0.3cm}\noindent
{\bf Theorem 4.5.} (Thm. 3.1 in \cite{Stu1}) {\it 
If a weight $\w$ defines a term order of the toric ideal $\IA$, then it
induces a regular triangulation $T_\w$. Moreover the Stanley-Reisner ideal
$SR_{T_\w}$ is equal to the radical of the initial ideal 
$\la LT_\w(\IA)\ra $. }

\vspace{0.3cm}
As an immediate corollary to this theorem, we see that the Gr\"obner fan
is a refinement of the fan $\{ {\cal C}({\cal A},T_\w) \}$. Since the cone
${\cal C}({\cal A},T_\w)$ decomposes according to (\ref{eqn:decomp}), 
we have similar
decomposition of ${\cal C}(\IA,\w)$ to ${\cal C}'(\IA,\w)$. In the
following we call the collection $\{{\cal C}'(\IA,\w)\}$ as the
Gr\"obner fan ${\cal F}(\IA)$. 

Now we determine a finite set of operators ${\cal D}_l$ which characterize the
power series $w_0(x,\rho;A)$ for each regular triangulation and a 
$\ZZ$-basis $A$ compatible with it. 
This provides us a way to analyze our resonant GKZ hypergeometric 
system. 

Let us consider a term order $\w$ of $\IA$. According to 
Th.4.5, the term oder $\w$ determines a regular triangulation
$T_\w$ and also a cone ${\cal C}'(\IA,\w) \subset {\cal C}'({\cal
A},T_\w)$. 
If the cone ${\cal C}'(\IA,\w)$ is simplicial and regular,
i.e., the integral generators of its one-dimensional boundary cones 
generate the lattice points ${\cal C}'(\IA,\w) \cap \Xi(\overline M)$, 
we simply make its dual cone ${\cal C}'(\IA,\w)^\vee$ and take the 
integral generators of this cone as a canonical $\ZZ$-basis $A$ of $L$ 
which is compatible with $T_\w$. If not, we
subdivide the cone ${\cal C}'(\IA,\w)$ into simplicial, regular cones
and reduce the problem to the former case. 
More generally, we may take a $\ZZ$-basis $A_\tau=\{
\l{1}_\tau,\cdots,\l{p-d}_\tau \}$ of $L$ compatible with $T_\w$
considering any simplicial, regular cone $\tau$ contained in ${\cal
C}'(\IA,\w)$ and making its dual $\tau^\vee$. 

Associated to $\w$, we have a Gr\"obner basis ${\cal B}_\w \subset \IA$. 
By B\"uchberger's algorithms to construct the (reduced) Gr\"obner basis, we
see that every generator $g \in {\cal B}_\w$ is a binomial of the form
$y^{l_+}-y^{l_-}$ with some $l \in L$. In the following, we assume 
${\cal B}_\w$ to be the reduced Gr\"obner basis which is determined 
uniquely for a term order $\w$ (, see Chapter 2 of \cite{CLO} for the 
properties of the reduced Gr\"obner basis). Translating this to differential
operator, we write the Gr\"obner basis ${\cal B}_\w=\{ {\cal D}_{l_1}, \cdots,
{\cal D}_{l_s} \} \; (1\leq s < \infty)$. Now, for each generator, we define 
\begin{equation}
J_l(\rho;A_\tau):= a_0 x_\tau^{-\rho} a^{l_\pm}\(\da{}\)^{l_\pm} 
            a_0^{-1}x_\tau^{\rho}
\;\;,
\label{eqn: Jl}
\end{equation}
where the choice in $l_\pm$ is made respectively by 
$\w\cdot l_+ - \w\cdot l_- >0 \;\; (<0) $.  
(The factor $a_0$ originate from the definition 
$w_0(x,\rho;A):= a_0 \Pi(a,\gamma)$ in Def.3.9.)

\vspace{0.3cm}\noindent
{\bf Definition 4.6.} 
For a term order $\w$ of $\IA$ and an arbitrary regular cone $\tau$ contained 
in ${\cal C}'(\IA,\w)$, we define, 
through the Gr\"obner basis ${\cal B}_\w=\{{\cal D}_{l_1},\cdots, {\cal D}_{l_s}\}$, 
an {\it indicial ideal} in $\CC[\rho_1,\cdots,\rho_{p-d}]$; 
\begin{equation}
Ind_\w(\tau)=\la  J_{l_1}(\rho,A_\tau),\cdots,J_{l_s}(\rho,A_\tau) \ra  \;\;. 
\end{equation}

Similarly to the indicial equations of the differential equations of 
Fuchs type, we also consider the {\it indicial equations} for our 
$\Ds$-hypergeometric system as algebraic equations for $\rho$ 
coming from the leading terms of the operators $D_l \;(l\in L)$. 
(Note that the leading term of an operator $D_l$ varies in general  
when a term order $\w$ varies. Here we consider for the indicial 
equations all possible leading terms when $\w$ varies inside $\tau$.)

\vspace{0.3cm}\noindent
{\bf Proposition 4.7.} {\it 
In the notation above, the indicial ideal $Ind_\w(\tau)$ 
coincides with the ideal generated by the indicial equations  
for the indices $\rho$ of the powerseries $w_0(x_\tau,\rho;A_\tau)$. }

\noindent
{\bf (Proof)} 
Consider an operator ${\cal D}_l \in {\cal B}_\w$. If 
$\w\cdot l_+ - \w\cdot l_- >0$, we multiply $a^{l_+}$ to obtain 
\begin{equation}
a^{l_+}{\cal D}_l =a^{l_+}\(\da{}\)^{l_+}
                   -a^{l_+-l_-}a^{l_-}\(\da{}\)^{l_-} \; .
\end{equation}
Since the initial ideal $\la LT_\w(\IA)\ra$ and thus 
the reduced Gr\"obner basis ${\cal B}_\w$ does not change for 
$\w$ in the interior of $\tau$, $Int(\tau)$, we have $\w\cdot(l_+-l_-)>0$ 
for all $\w\in Int(\tau)$, i.e., $l_+-l_-\in \tau^\vee\cap L$. 
Since we have chosen the $\ZZ$-basis 
$A_\tau=\{ \l{1}_\tau,\cdots,\l{p-d}_\tau\}$ so that it generates 
all integral points in $\tau^\vee \cap L$,  
$a^{l_+-l_-}$ is a monomial of $x_\tau$, which
vanish in the limit $x_\tau \rightarrow 0$. The same argument applies 
to the case $\w\cdot l_+ - \w\cdot l_- <0$. Therefore the indicial equations 
arising from the operators $D_l \in {\cal B}_\w$ exactly coincide 
with the generators of the indicial ideal (\ref{eqn: Jl}). 
For general operators $D_l \; (l\in L)$ , depending on 
the weight $\w \in Int(\tau)$, we have two possible leading terms. 
However for both of them, 
owing to the defining property of the Gr\"obner basis, 
we have $LT_\w(D_l)=\(\da{}\)^\mu LT_\w(D_{l_k})$
for some $k$ and $\mu$. Multiplying a monomial $a^{\mu+l_{k\pm}}$, we
obtain 
\begin{equation}
a_0x_\tau^{-\rho} a^{\mu+l_{k\pm}} LT_\w(D_l) a_0^{-1} x_\tau^{\rho} 
=F(\rho) \, J_{l_k}(\rho;A_\tau) \;\;,
\end{equation}
with some polynomial $F(\rho)$.  
Thus we see all polynomial relations of $\rho$ comming from the 
leading terms are in $Ind_\w(\tau)$. 
\par
Conversely, since all generators of the ideal $Ind_\w(\tau)$  
give the indicial equations related to ${\cal B}_\w$, 
the ideal $Ind_\w(\tau)$ is contained in the other. 
Therefore the two ideals are the same.  \hfill $\Box$

\vspace{0.3cm}
Now based on Prop.4.7, we may claim the following; 

\vspace{0.3cm}\noindent
{\bf Proposition 4.8.} {\it 
Consider a compact toric variety $\PP_{{\cal F}(\IA)}$ associated to 
the Gr\"obner fan $({\cal F}(\IA), \Xi(\overline M))$. Then
for any resolution $\PP_{\tilde{\cal F}(\IA)} \rightarrow \PP_{{\cal
F}(\IA)}$ associated to a refinement $(\tilde{\cal F}(\IA), \Xi(\overline M)) 
\rightarrow ({\cal F}(\IA), \Xi(\overline M))$, we have integral powerseries 
of the form $w_0(x_\tau,\rho;A_\tau)$ $(\rho\in V(Ind_\w(\tau)))$ 
at each boundary point given by the normal
crossing toric divisors, namely at the origin of 
$Hom_{s.g.}(\tau^\vee\cap L,\CC)$. We will call this compactification
Gr\"obner compactification. }

\vspace{0.3cm}\noindent
{\bf Remark.} 
Since Prop.4.7 provides us only a necessary condition for the indices 
$\rho$ to give a powerseries solution $w_0(x,\rho;A_\tau)$, we do not claim 
by Prop.4.8, although we expect, that all $\rho\in V(Ind_\w(\tau))$ 
form the powerseries solutions of our $\Ds$-hypergeometric system.

\newsubsec{(4-3) Resonance of $\Ds$-hypergeometric system} 
When the polytope $\Ds$ is of type I or II, we have seen in 
the Remark right after Th.3.8 that the 
$\Ds$-hypergeometric system becomes ``maximally resonant'' for a maximal 
triangulation $T_o$. Here we study this resonance in detail restricting 
our attention to the polytopes of type I or II. We also comment about 
the case of type III. 

We call a collection of vertices ${\cal P}=\{\bn{i_1},\cdots,\bn{i_a}\}$ 
{\it primitive} if ${\cal P}$ does not form a simplex in $T_o$ but 
${\cal P}\setminus \{ \bn{i_s}\}$ does for any $\bn{i_s}\in {\cal P}$. 
By definition of the Stanley-Reisner ideal, it is easy to deduce that 
the monomials that corresponds to primitive collections generate the 
ideal $SR_{T_o}$. 

Let us denote by $\Sigma(1\times\Ds,T_o)$ the fan in $\overline N_\RR$ that 
is naturally associated to the triangulation $T_o$ of $P$. Since the 
volumes of all $d+1$ simplices in $T_o$ are unimodular for the polytope 
$\Ds$ of type I or II, the fan $\Sigma(1\times\Ds,T_o)$ consists of regular 
cones. Therefore if we have a primitive collection 
${\cal P}=\{\bn{i_1},\cdots,\bn{i_a}\}$, we obtain
\begin{equation}
\bn{i_1}+\cdots+\bn{i_a}=\sum_k c_k \bn{j_k} \;\; (c_k \in \ZZ_{\geq0})
\label{eqn: prim}
\end{equation}
where $\{\bn{j_k} \vert c_k\not= 0\}$ generates a cone that contains 
the vector in the left hand side. 
Writing (\ref{eqn: prim}) as 
$\bn{i_1}+\cdots+\bn{i_a}-\sum c_k \bn{j_k}=0$, we read 
the corresponding {\it primitive relation} $l({\cal P}) \in L$.

\vspace{0.3cm}\noindent
{\bf Lemma 4.9.} {\it 
Every primitive collection of a maximal triangulation $T_o$ does not 
contain the point $\bn{0}=1\times \vec 0$. }
\par
\noindent
{\bf (Proof)} 
Suppose a primitive collection is given by ${\cal P}=\{ \bn{0}, \bn{i_1},
\cdots,\bn{i_a} \} $ $(1\leq i_1,\cdots,i_a \leq p)$. Since it is primitive,  
the simplex $\la  \bn{i_1},\cdots,\bn{i_a} \ra $ must be a simplex in  
the triangulation $T_o$, which means that this simplex is a face of 
some maximal dimensional simplex in $T_o$. Since $T_o$ is a maximal 
triangulation in which every maximal dimensional simplex contains 
the vertex $\bn{0}$, we see the simplex 
$\la \bn{0},\bn{i_1},\cdots,\bn{i_a} \ra$ must be a simplex 
in $T_o$, which is a contradiction. \hfill $\Box$

\vspace{0.3cm}\noindent
{\bf Proposition 4.10.} {\it 
For a term order 
$\w$ such that $T_\w$ is a maximal triangulation, 
the initial ideal $\la LT_\w(\IA)\ra $ is radical and 
$\la LT_\w(\IA)\ra =SR_{T_\w}$. }
\par\noindent
{\bf (Proof)} 
Consider the primitive collections for the triangulation $T_\w$,  
which generate the Stanley-Reisner ideal $SR_{T_\w}$. Write a primitive 
collection ${\cal P}=\{ \bn{i_1},\cdots,\bn{i_a}\}$ and the corresponding 
primitive relation as $l({\cal P})$ considering the relation 
(\ref{eqn: prim}).  
For a term order $\w$, the regular triangulation $T_\w$ 
is induced from the lower faces of the polytope 
$P_\w={\rm Conv.}(\,\{\tilde\ns{0}, \cdots ,\tilde\ns{p} 
\;\vert\; \tilde\ns{k}=$ $\w_k\times\ns{k}\;(k=0,\cdots,p)\}\,)$. 
Then the convex hull 
${\rm Conv.}(\{\tilde\ns{i}\;\vert\; \bn{i}\in{\cal P}\})$ is not a simplex  
that corresponds to a lower face of $P_\w$. Therefore 
we have a ``height'' inequality 
$( \tilde\ns{i_1}+\cdots+\tilde\ns{i_a} )_1 > (\sum c_k \tilde\ns{j_k})_1$, 
namely,
\begin{equation}
\w_{i_1}+\cdots+\w_{i_a} > \sum_k c_k \w_{j_k} \;\;.
\label{eqn: height}
\end{equation}
This means that $LT_\w(y^{l({\cal P})_+}-y^{l({\cal P})_-})=
y_{i_1}\cdots y_{i_a}$, which is one of the generators of the ideal 
$SR_{T_\w}$. Since this argument applies to all primitive collections, 
we conclude $SR_{T_\w}\subset \la LT_\w(\IA)\ra $. Since the opposite 
inclusion follows from $SR_{T_\w}=\sqrt{ \la LT_w(\IA)\ra  }$ (Th.4.5), 
we conclude $SR_{T_\w}=\la LT_w(\IA)\ra $, which proves the initial 
ideal is radical. \hfill $\Box$

\vspace{0.3cm}\noindent
{\bf Corollary 4.11.} {\it 
Under the hypothesis in the previous proposition, 
the set of all possible primitive collections 
$\{ {\cal P}_1,\cdots,{\cal P}_s \}$ of $T_\w$ determines the 
Gr\"obner basis by ${\cal B}_\w=\{ {\cal D}_{l({\cal P}_1)}, \cdots, 
{\cal D}_{l({\cal P}_s)} \}$. 
And the indicial ideal $Ind_\w(\tau)$ 
is homogeneous for an arbitrary regular cone $\tau$ contained 
in ${\cal C}'(\IA,\w)$. }
\par\noindent
{\bf (Proof)} 
By definition, $SR_{T_\w}$ is generated by the monomials 
corresponding to primitive collections. From the argument in the 
proof of Prop.4.10, we know 
$SR_{T_\w}=\la LT_\w( {\cal D}_{l({\cal P}_1)}),$ 
$\cdots, LT_\w( {\cal D}_{l({\cal P}_s)}) \ra$. 
This combined with $\la LT_\w(\IA) \ra =SR_{T_\w}$ establishes 
that ${\cal B}_\w$ is the Gr\"obner basis.   
{}For the rest, we note that any primitive collection ${\cal P}=
\{ \bn{i_1},\cdots,\bn{i_a}\}$ does not contain $\bn0$ according to 
Lemma 4.9. Now we have 
\begin{equation}
J_{l({\cal P})}(\rho;A_\tau)=a_0x_\tau^{-\rho} 
a^{l({\cal P})_+} \(\da{}\)^{l({\cal P})_+} a_0^{-1}x_\tau^\rho 
= x_\tau^{-\rho} \ta{{i_1}}\cdots\ta{{i_a}} x_\tau^\rho ,
\end{equation}
where $ \ta{{}} =a \da{} $, which proves that the generator 
$J_{l({\cal P})}(\rho;A_\tau)$ is 
homogeneous in $\rho$.   \hfill $\Box$

\vspace{0.3cm}\noindent
{\bf Remark.} 
If we combine a general result that the GKZ system is holonomic 
\cite{GKZ1}, i.e., its solution space is finite dimensional, 
with our Corol. 4.11, we may 
conclude that the zero is the only solution for the indices $\rho$. 
This is the maximal $T$-resonance in our approach.  We will give 
following \cite{HLY2} an independent proof about this in the next section. 

As we remarked before, our $\Ds$-hypergeometric 
system for the polytope $\Ds$ of 
type III does not share this property. Here we can explain the difference.  
We first note that the primitive collections 
generate the Stanley-Reisner ideal and has the property in Lemma 4.9 
irrespective to the type of polytopes. 
The only change in the above arguments 
is in the definition of the primitive relation.  Namely, since all  
cones are not regular in type III case, for some primitive collection 
the equation (\ref{eqn: prim}) should be replaced by 
\begin{equation}
\lambda_{i_1}\bn{i_1}+ \cdots + \lambda_{i_a} \bn{i_a} 
= \sum_k c_k \bn{j_k} \;\;,
\label{eqn: sing}
\end{equation}
with some positive integers $\lambda_{i_1},\cdots,\lambda_{i_a}$ not all 
equal to one.  
Accordingly the leading term $LT_\w(y^{l({\cal P})_+}-y^{l({\cal P})_-})$ 
will be replaced by $(y_{i_1})^{\lambda_{i_1}}\cdots 
(y_{i_a})^{\lambda_{i_a}}$. This indicates that the initial ideal 
$LT_\w(\IA)$ is no longer radical, and therefore the generators 
$J_l(\rho;A_\tau)$ become inhomogeneous. When translating the monomial 
$y^{l({\cal P})_+}$ to the differential operator 
$a^{l({\cal P})_+} \( \da{} \)^{l({\cal P})_+}$, each $\lambda_i$-fold 
degeneration to zero 'splits' to simple zeros. 
Thus every index does not degenerate to zero, 
although we still have a simple zero.

\section{5. Large compex Structure Limit}

Here we will study in detail the maximal resonance of the 
$\Ds$-hypergeo-metric system. We will identify this resonance with the
large complex structure limit (LCSL), i.e., a celebrated boundary point
in the moduli space of Calabi-Yau manifolds\cite{Mor}. 

\newsubsec{(5-1) Maximal degeneration} In this subsection, we will restrict
our arguments to the polytopes of type I or II. In these two cases, 
we have a nonsingular projective toric variety
$\PP_{\Sigma( \Ds,T_o)}$ for a maximal triangulation. We focus on
the Chow ring of this toric variety. The Chow ring of a variety is a
free abelian group generated by irreducible closed subvarieties, modulo
rational equivalence, which is endowed with the ring structure via the
intersection products. In case of non-singular compact toric varieties
$\PP_{\Sigma(\Ds,T_o)}$, it has a simple description in terms of the
(toric) divisors;

\vspace{0.3cm}\noindent
{\bf Proposition 5.1.} (sect.3.3 of \cite{Oda}, sect.5.2 of \cite{Ful}) 
{\it 
The Chow ring $A^*(\PP_{\Sigma(\Ds,T_o)})$ is isomorphic to the
cohomology ring $H^{2*}(\PP_{\Sigma(\Ds,T_o)},\ZZ)$ and is given by
\begin{equation}
A^*(\PP_{\Sigma(\Ds,T_o)})=\ZZ[D_0,\cdots,D_p]/(SR_{T_o}+\bar R), 
\end{equation}
where $D_k \; (k>0)$ represents the toric-divisor determined by 
the integral point $\ns{k}$. $SR_{T_o}$ is the Stanley-Reisner 
ideal and $\bar R$ is the ideal generated by linear relations 
$\sum_{k=0}^p \la 1\times u, \bn{k}\ra  D_{k} =0 \;\;(u\in M)\;. $ }

\vspace{0.3cm}\noindent
{\bf Note.} 
Owing to lemma 4.9, we can take the generators of $SR_{T_o}$ that 
do not contain $D_0$. Therefore the generator $D_0$ plays only a dummy 
role, although it makes sense as a divisor if we consider a toric 
variety defined by the (non-complete) fan 
$\Sigma(1\times\Ds,T_o) \subset \overline N_\RR$. 
 
\vspace{0.3cm}
Now consider a term order $\w$ of the toric ideal $\IA$ and denote the 
Gr\"obner basis by ${\cal B}_\w=\{ {\cal D}_{l_1},\cdots,{\cal D}_{l_s}\}$. We define
\begin{equation}
I_l(\ta{{}}):=a_0 a^{l_\pm}\(\da{{}} \)^{l_\pm} a_0^{-1} 
\end{equation}
for each $LT_\w({\cal D}_l)=(\da{{}})^{l_\pm}$ in a similar 
way to $J_l(\rho;A_\tau)$. Obviously these two are related by
$J_l(\rho;A_\tau)=x_\tau^{-\rho}I_l(\ta{{}} )x_\tau^{\rho}$. 
We consider the following ideals in $\CC[\ta0,\cdots,\ta{p}]$,
\begin{equation}
I_\w:=\la  I_{l_1}(\ta{{}} ),\cdots, I_{l_s}(\ta{{}} ) \ra  \;,\; 
\bar R_a:=\la \sum_{i=0}^p\la  1\times u, \bn{i}\ra \ta{i}
          \;\vert\; u\in M \ra  \;\;.
\end{equation}

\vspace{0.3cm}\noindent
{\bf Proposition 5.2.} {\it 
For a term order $\w$ of $\IA$ and an arbitrary regular cone $\tau$
contained in ${\cal C}'(\IA,\w)$, we have 
\begin{equation}
\CC[\rho]/Ind_\w(\tau) \cong \CC[\ta{{}}]/(I_\w+\bar R_a) \;\;.
\end{equation} }
\par
\noindent
{\bf (Proof)} 
When we take the $\ZZ$-basis $A_\tau=\{ \l1_\tau, \cdots,\l{p-d}_\tau
\}$, we have $\ta{i}=\sum_{k=1}^{p-d} (\l{k}_\tau)_i
\theta_{x^{(k)}_\tau}$. Then the homomorphism $\phi : \CC[\ta{{}}]\rightarrow 
\CC[\theta_{x_\tau}]\cong\CC[\rho]$ induced by this relation is surjective, 
since rank$(L)=p-d$, and satisfies 
$\Ker \phi = \bar R_a $ and $\phi(I_\w)= Ind_\w (\tau)$. 
This proves the assertion. 
\hfill $\Box$

\vspace{0.3cm}\noindent
{\bf Proposition 5.3.} {\it 
Consider a term order $\w$ of $\IA$ with $T_\w$ a maximal
triangulation $T_o$. Then for any regular cone $\tau$ contained in 
${\cal C}'(\IA,\w)$, the variety associated to the indicial ideal 
$Ind_\w(\tau)$ consists only one point, i.e., 
\begin{equation}
V(Ind_\w(\tau))=\{0\} 
\;\;. 
\end{equation}
  }
\par
\noindent
{\bf (Proof)} 
By Corol. 4.11, we know the indicial ideal is homogeneous 
for a term order $\w$ of the given property. Moreover the 
ideal $I_\w$ coincides with the Stanley-Reisner ideal $SR_{T_o}$. 
Therefore we have
\begin{equation}
\CC[\rho]/Ind_\w(\tau) \cong \CC[\ta{{}}]/(I_\w+\bar R_a) 
\cong A^*(\PP_{\Sigma(\Ds,T_o)})\otimes\CC \;\;, 
\end{equation}
which is finite dimensional. Since the ideal is homogeneous, 
the claim follows.    \hfill $\Box$

We write our series (\ref{eqn: wnot}) for generic $\rho$ by $w_0(x,\rho;A)=
\sum_{m\in \ZZ^{p-n}_{\geq 0}} c(m+\rho)x^{m+\rho}$. As remarked after 
Def.3.9, the value $\rho=0$ might cause the infinity in the numerator for 
some coefficient $c(m+\rho)$. One way to treat this infinity problem is 
to take the limit $\rho\rightarrow 0$ as discussed there and 
we will come back to this recepie in Claim 5.8.  Here following 
\cite{HLY2}, we introduce the series
\begin{equation}
w_0(x_\tau,0;A_\tau)_{\geq0}:= 
\sum_{m\in\ZZ^{p-n}_{\geq0},\;  -\sum_k m_k\l{k}_\tau\geq0} c(m) x^m \;\;.
\label{eqn: wnotpositive}
\end{equation}

\vspace{0.3cm}\noindent
{\bf Theorem 5.4.} (Th.5.2 in \cite{HLY2}) {\it 
For a term order $\w$ with $T_\w$ a maximal triangulation 
and any regular cone $\tau$ contained in ${\cal C}'(\IA,\w)$, 
the series $w_0(x_\tau,0; A_\tau)_{\geq0}$ is
the only powerseries solution of the $\Ds$-hypergeometric system 
about the origin of $U_\tau=Hom_{s.g.}(\tau^\vee\cap L,\CC)$. }

\vspace{0.3cm}
To prove this theorem, we prepare the following lemma;

\vspace{0.3cm}\noindent
{\bf Lemma 5.5.} {\it 
Consider a subset $S\not=\{\phi\}$ that is contained in $\tau^\vee\cap L$.
There exist an element $\delta\in S$ and a simplicial, regular cone 
$C_\delta^\vee \subset L_\RR$ such that $C_\delta^\vee$ contains
both the subset $S-\delta$ and the cone $\tau^\vee$. }
\par\noindent
{\bf (Proof)} 
Consider a hyperplane $H(v;z_0)$ with a normal vector $v\in \tau$ and
passing through a point $z_0$ in $\tau^\vee$. When we consider a parallel
transport $H(v,t z_0)$ $(t\geq0)$ of the hyperplane, we may find  
the minimal $t_0$ such that $H(v, t_0 z_0)\cap S \not=\{\phi\}$ while  
$H(v,t z_0)\cap S =\{\phi\} \;\; (t <t_0)$. Changing the normal vector
$v$ slightly, if necessary, we may assume the intersection 
$H(v, t_0 z_0)\cap S$ occurs at a point $\delta$. Now for this $\delta$,
we see that the union $U:=(\tau^\vee\cap L) 
\cup (S-\delta) \setminus\{0\}$ is contained in the half space $H_>(v,0)$.  
Therefore the normal cone to the set $U$ at the origin is 
strongly convex, polyhedral cone. Since a strongly convex, polyhedral cone 
can be inside a simplicial, regular cone, the assertion follows. \hfill $\Box$

\noindent
{\bf (Proof of Th.5.4.)}
To prove the theorem, we write the series $w_0(x_\tau,0;$ $A_\tau)$ in
terms of $a_0,\cdots,a_p$ by
\begin{equation}
w_0(a,0,\tau)=\sum_{l \in \tau^\vee\cap L} c_l a^l \;\; ,
\end{equation}
with $c_0=1$. Now suppose we have two different series of this form. 
Then the difference of the two may be written by  
$r(a,0,S)=\sum_{l \in S} d_l a^l $ 
with a subset $S\subset \tau^\vee \cap L\setminus \{0\}$. Using the
result in the lemma 5.5, we may write this series via nonzero $\delta$ as
\begin{equation}
r(a,0,C_\delta)=a^\delta \sum_{l \in C_\delta^\vee \cap L} 
                d_{l+\delta}a^{l} \;\;,
\label{eqn: delta}
\end{equation}
or $r(x_\tau,0,A_{C_\delta})=x_\tau^{\rho(\delta)} \sum_{n \in
\ZZ^{p-d}_{\geq 0}} d(n) x_\tau^n $ with $\rho(\delta)\not=0 ,d(0)\not=0$ and 
$C_\delta \subset \tau$. This is a contradiction to Prop.5.3. \hfill $\Box$

\vspace{0.3cm}\noindent
{\bf Remark.} 
By direct evaluation of the period integral (\ref{eqn: Pi})\cite{Bat2}, 
we can verify that $a_0\Pi(a)$ exactly coincides with the powerseries in
Th.5.4 when expressed in terms of the $\ZZ$-basis $A_\tau$ (Prop.5.15 
\cite{HLY2}).

\newsubsec{(5-2) All solutions about maximal degeneration points} 
Here we determine other solutions about maximal
degeneration points, all of which contains logarithmic singularities. 
As in the previous subsection, our arguments are restricted to the
polytopes of type I or II. 

Let us note that the first degree elements of the Chow ring,
$A^1(\PP_{\Sigma(\Ds,T_o)})$, describe the Picard group of the toric
variety $\PP_{\Sigma(\Ds,T_o)}$ and may be expressed by
\begin{equation}
A^1(\PP_{\Sigma(\Ds,T_o)})=\ZZ D_0\oplus\cdots\oplus\ZZ D_p / \bar R 
\;\cong\; \Xi(\overline M) \;\;.
\end{equation}
From this we see a dual pairing between the Picard group and the lattice
$L\cong \Xi(\overline N)$;
\begin{equation}
A^1(\PP_{\Sigma(\Ds,T_o)})\times L \rightarrow \ZZ \;\;.
\label{eqn: dualpL}
\end{equation}

\vspace{0.3cm}\noindent
{\bf Definition 5.6.} 
For a $\ZZ$-basis $A_\tau=\{ \l{1}_\tau, \cdots, \l{p-d}_\tau \}$ of $L$
determined from a term order $\w$ with $T_\w$ equal to a maximal 
triangulation $T_o$, we denote its dual by 
$A_\tau^\vee=\{ J_{\tau,1},\cdots,J_{\tau,p-d} \}$
or simply by $\{ J_1\cdots,J_{p-d}\}$ when its dependence on $\tau$ is
obvious. 

\vspace{0.3cm}\noindent
{\bf Note.} 
By construction, the basis $A_\tau^\vee$ consists of 
the integral generators of the simplicial, 
regular cone $\tau$ contained in ${\cal C}'(\IA,\w)$ ($= 
{\cal C}'({\cal A},T_o)$ by Prop.4.10). ${\cal C}'({\cal A},T_o)$ consists 
of convex functions on $T_o$ which may be identified with the convex 
functions on the fan $\Sigma(\Ds,T_o)$. Since the set of all convex functions 
on the fan $\Sigma(\Ds,T_o)$ determines the closure of the K\"ahler cone 
of $\PP_{\Sigma(\Ds,T_o)}$ (see Corol.2.15 in \cite{Oda}), 
the bases $J_{\tau,1},\cdots, J_{\tau,p-d}$ generate a simplicial, 
regular cone contained in this closure of the K\"ahler cone.

\vspace{0.3cm}{\noindent}
{\bf Definition 5.7.} 
For the powerseries $w_0(x_\tau,\rho;A_\tau)=
\sum_{n\in\ZZ_{\geq0}^{p-d}} c(n+\rho) 
x_\tau^{n+\rho}$ in (\ref{eqn: wnot}), we
define 
\begin{equation}
w_0(x_\tau,J;A_\tau):= \sum_{n\in\ZZ_{\geq0}^{p-d}} 
c\(n+{J \over 2\pi i}\) x_\tau^{n+{J\over 2\pi i}} \;\;, 
\label{eqn: wJ}
\end{equation}
as the Taylor series expansion of $w_0(x_\tau,\rho;A_\tau)$ about $\rho=0$ 
followed by the substitution $\rho={J\over 2\pi i}$, where $J$'s 
are defined in Def.5.6.

\vspace{0.3cm}
In refs.\cite{HKTY1}\cite{HKTY2}\cite{HLY1}\cite{HLY2}, 
it is widely verified

\vspace{0.3cm}\noindent
{\bf Claim 5.8.} {\it  
The expansion (\ref{eqn: wJ}) exists as an element in 
$A^*(\PP_{\Sigma(\Ds,T_o)})\otimes \CC\{x_\tau\}[\log x_\tau]$, 
and the coefficient series constitute a complete set of 
the local solutions about the maximal degeneration points. 
Especially the limit $w_0(x_\tau,\rho;A_\tau)\vert_{\rho\rightarrow0}$ 
coincides with $w_0(x_\tau,0; A_\tau)_{\geq0}$ in Th.5.4. }

\vspace{0.3cm}\noindent
{\bf Remark.} 
We comment about the case of the polytopes of type III. 
In this case, since the toric variety $\PP_{\Sigma(\Ds,T_o)}$ is singular, 
the Chow ring should be considered over $\QQ$. Under this modification the
expansion (\ref{eqn: wJ}) makes sense in 
$A^*(\PP_{\Sigma(\Ds,T_o)})_\QQ \otimes \CC\{x_\tau\}[\log x_\tau]$. 
Then the coefficient series should be in a subspace of
the whole solution space of the $\Ds$-hypergeometric system. More
precisely, as we see in Remark after Corol. 4.11, the initial ideal 
$LT_\w(\IA)$ is no longer radical but we have a strict inclusion 
$LT_\w(\IA) \subset \sqrt{LT_\w(\IA)}$. 
As discussed there, we have 
$LT_\w({\cal D}_l)=(\da{i_1})^{\lambda_{i_1}} 
\cdots (\da{i_a})^{\lambda_{i_a}}$ 
for some element of the Gr\"obner basis ${\cal B}_\w=\{ {\cal D}_{l_1},
\cdots, {\cal D}_{l_s} \}$. If we define ${\rm rad}(LT_\w({\cal D}_l))
:=\da{i_1}\cdots\da{i_a}$, 
then the radical may be expressed by 
$\sqrt{LT_\w(\IA)}=\la  {\rm rad}(LT_w({\cal D}_{l_1})), 
\cdots, {\rm rad}(LT_w({\cal D}_{l_s})) \ra $. Correspondingly, if we define 
$\tilde I_\w(\ta{{}}):=
a_0 a_{i_1}\cdots a_{i_a} {\rm rad}(LT_\w({\cal D}_l)) a_0^{-1}$, 
we naturally come to the ``radical'' of the indicial ideal 
$\tilde{Ind}_\w(\tau):=\la  \tilde J_{l_1}(\rho;A_\tau), \cdots, \tilde 
J_{l_s}(\rho;A_\tau) \ra $ with $\tilde J_l(\rho;A_\tau):=x_\tau^{-\rho}\tilde 
I_\w(\ta{{}}) x_\tau^{\rho}$. By definition, we have strict inclusions 
$Ind_\w(\tau)$ $\subset \tilde{Ind}_\w(\tau)$ and $V(Ind_\w(\tau)) \supset 
V(\tilde{Ind}_\w(\tau))$. As is clear now, our Prop.5.2 and Prop.5.3 
apply to the ``radical'' $\tilde{Ind}_\w(\tau)$ under the replacements 
$I_\w$ by $\tilde I_\w$ and the Chow ring by that over $\QQ$. 
The expansion (\ref{eqn: wJ}) gives all logarithmic solutions 
which arise from the degeneration $V(\tilde{Ind}_\w(\tau))$ $=\{ 0\}$.

\newsubsec{(5-3) LCSL of Calabi-Yau hypersurfaces} 
So far we have been concerned with the $\Ds$-hypergeometric system. 
Since the period integral 
(\ref{eqn: Pi}) of Calabi-Yau hypersurface $X_\Ds$ satisfies the 
(extended) $\Ds$-hypergeometric system, a complete set of 
the period integrals of $X_\Ds$ 
should be found in the set of solutions of the $\Ds$-hypergeometric system. 
We will find that the expansion (\ref{eqn: wJ}) contains the period integrals 
in a natural way from the mirror symmetry. 

Before going into this topic, we need to discuss about the compactification 
of the moduli space ${\cal M}(X_\Ds(a))$ of the polynomial deformation 
of the Calabi-Yau hypersurface $X_\Ds$. Through a detailed analysis of 
the local solutions of the $\Ds$-hypergeometric system, we have arrived 
at a natural compactification, the Gr\"obner compactification 
$\PP_{{\cal F}(\IA)}$ in Prop.4.8. Now it is natural to adopt this 
compactification as that of the moduli space ${\cal M}(X_\Ds(a))$. 
However, one problem arises when the hypersurface (, precisely its ambient 
space,) has non-trivial automorphisms. We need to mod out the space 
$\PP_{{\cal F}(\IA)}$ by the induced actions from the automorphisms, 
whose infinitesimal forms are described in (\ref{eqn: extGKZ}). 
Here to avoid getting involved in the problems related to the 
actions of the automorphisms, we take a ``gauge choice'' that sets to 
zero all polynomial deformations corresponding to integral points 
on codimension-one faces of $\Ds$. Note that, in view of Prop.2.6,  
the degree of the freedom associated to the non-trivial automorphisms 
would be fixed by this gauge choice. 
In the following, we use the subscript $s$ 
(, $s$ of simply-minded!,) to indicate this 
naive choice of the ``gauge''; for example $\Delta_s^*$-hypergeometric system, 
the toric ideal $\IAs$ etc.  Note that all the polytopes of type II 
will be treated as the polytopes of type III under this prescription. 

\vspace{0.3cm}\noindent
{\bf Definition 5.9.}
As an compactification of ${\cal M}(X_\Ds(a))$, we define 
\begin{equation}
\overline{\cal M}(X_{\Ds}(a))=\PP_{{\cal F}(\IAs)} \;\;.
\end{equation}

\vspace{0.3cm}
Now we consider the toric part of the Chow ring of the Calabi-Yau 
hypersurface $X_{\D_s}$, which comes from the ambient space by 
restriction. Since we have $[X_{\D_s}]=D_1+\cdots+D_p$ for the divisor of 
the Calabi-Yau hypersurface, the restriction may be attained by the 
quotient as follow; 

\vspace{0.3cm}\noindent
{\bf Definition 5.10.}
\begin{equation}
A^*(X_{\D_s})_{toric}=
A^*(\PP_{\Sigma(\Delta_s^*,T_o)})_\QQ/ Ann(D_1+\cdots+D_p) \;\;, 
\label{eqn: chowX}
\end{equation}
where $Ann(x)$ is defined by 
$Ann(x)=\{ \; y\in {\cal R} \;\vert\; x\,y=0 \;\}$ for a ring ${\cal R}$.

\vspace{0.3cm}\noindent
{\bf Claim 5.11.} {\it 
Period integrals about a LCSL of Calabi-Yau hypersurface $X_{\Delta_s^*}$ are 
extracted from the series $w_0(x_\tau,J;A_\tau)$ (\ref{eqn: wJ}) 
expanded in $A^*(X_{\D_s})_{toric}\otimes \CC\{x_\tau\}[\log x_\tau]$. }

\vspace{0.3cm}\noindent
{\bf Remark.}
In general, the period integrals of Calabi-Yau hypersurfaces satisfy the 
differential equations of Fuchs type, so-called the Picard-Fuchs equations
\cite{PF}. 
Picard-Fuchs equations determines the period integrals as its solutions. 
Our Claim 5.11 says that our $\Ds$-hypergeometric system is reducible in 
general and contains the Picard-Fuchs equation as a component of it. 
In refs.\cite{HKTY1}\cite{HLY1}, it is verified in several examples 
of all types of the polytopes that Picard-Fuchs 
equations are derived from the (extended) $\Ds$-hypergeometric system 
after a factorization of the operator $\ta{1}+\cdots+\ta{p}$ from the 
left, which we may identify with the quotient by $Ann(D_1+\cdots+D_p)$ in 
(\ref{eqn: chowX}).

\section{6. Prepotential}

In this section, we study so-called the 
prepotential \cite{Saito} near a  LCSL in detail. 
Under the mirror map, a LCSL is mapped to a large 
radius limit in which the instanton corrections to the prepotential 
are suppressed exponentially, and has important applications to 
the enumerative geometry. Also the prepotential determines the 
{\it special K\"ahler geometry} on the moduli space 
$\overline {\cal M}(X_{\Delta_s^*})$ and  
that on the complexified K\"ahler moduli space of the 
mirror $X_{\D_s}$. 

In this section, we fix a term order $\w$ for which $T_\w$ is a maximal 
triangulation of $\Delta_s^*$ and take the $\ZZ$-basis $A_\tau$ choosing a 
regular cone $\tau$ in ${\cal C}'(\IAs,\w)$. 
Based on Claim 5.11, 
we expand the series $w_0(x_\tau,J;A_\tau)$ defined in (\ref{eqn: wJ}) 
(see also \cite{Sti}) as; 
\begin{equation}
w_0(x_\tau,J;A_\tau)=
w^{(0)}(x_\tau,J)+w^{(1)}(x_\tau,J)+{1\over2!}w^{(2)}(x_\tau,J)
+{1\over3!}w^{(3)}(x_\tau,J) \;,
\label{eqn: wJexp}
\end{equation}
where the superscripts indicate the degree in the Chow ring 
$A^*(X_{\D_s})_{toric}\otimes \CC\{x_\tau\}[\log x_\tau]$. 

\vspace{0.3cm}\noindent
{\bf Definition 6.1.} 
The {\it special coordinate} $(t_1,\cdots, t_{p-d})$ of 
the special K\"ahler geometry is defined 
by the ratios of the period integrals;
\begin{equation}
t\cdt J = {w^{(1)}(x_\tau,J) \over w^{(0)}(x_\tau,J) } \;\;,
\label{eqn: special}
\end{equation}
where we abuse the letters $J_1, \cdots, J_{p-d}$ to represent the images 
of the $J$'s under the quotient (\ref{eqn: chowX}). 
The inverse series of this relation will be called the {\it mirror map}. 

\vspace{0.3cm}\noindent
{\bf Note.} Since $w^{(1)}(x_\tau,J)$ is linear in log$x_\tau$, the mirror 
map takes the form $x_\tau(q):=x_\tau(q_1,\cdots,q_{p-d})$ with 
$q_k:={\rm e}^{2\pi i t_k}$. It is easy to see that $x_\tau^{(k)}(q)
=q_k(1+{\cal O}(q))$. 

\vspace{0.3cm}\noindent
{\bf Definition 6.2.} 
We define the prepotential in the special coordinate by 
\begin{equation} 
F(t):=\int_{X_{\D_s}} {\cal F}(x_\tau(q), J) \;\;,
\label{eqn: Ft}
\end{equation}
with the {\it invariant} density 
\begin{equation}
{\cal F}(x_\tau(q), J)={1\over2} \({1\over w^{(0)}}\)^2 
\{w^{(0)}(-{1\over3!}w^{(3)}-{c_2(X_{\D_s})\over 12}w^{(1)})+
w^{(1)}({1\over2!}w^{(2)}) \} \;.
\label{eqn: Fd}
\end{equation}
The integration symbol 
$\int_{X_{\D_s}} :=\int_{\PP_{\Sigma({\Delta_s^*},T_o)}} 
[X_{\D_s}]$ is meant to take the coefficient of the 'volume form' in the 
Chow ring $A^*(\PP_{\Sigma(\Delta_s^*,T_o)})_\QQ$ normalized by 
$\int_{\PP_{\Sigma({\Delta_s^*},T_o)}} [X_{\D_s}]c(\PP_{\Sigma(\Delta_s^*,T_o)})/(1+[X_{\D_s}])=\chi(X_{\D_s})$. (For the normalization when 
$\chi(X_{\D_s})=0$, see ref.\cite{HLY1}.) 

\vspace{0.3cm}\noindent
{\bf Note.} 
It would be instructive to summarize the general description\cite{Str} of 
the special K\"ahler geometry on the complex structure moduli space of 
Calabi-Yau threefolds. Let us denote the holomorphic 3-form of 
a family of Calabi-Yau threefolds $W_\psi$ by $\Omega(\psi)$. 
We take a symplectic basis $\{A_a,B_b\}$  $(a,b=0,\cdots,h^{2,1}(W))$ 
of $H_3(W,\ZZ)$ and construct the period integrals $z_a(\psi)=\int_{A_a}
\Omega(\psi)$ and ${\cal G}_b(\psi)=\int_{B_b}\Omega(\psi)$. 
Then the holomorphic 3-form may be written by $\Omega(\psi)=\sum_a z_a(\psi)
\alpha_a + \sum_b {\cal G}_b(\psi)\beta_b$ in terms of the dual 
bases $\alpha_a$ and $\beta_b$ in $H^3(W,\ZZ)$. Locally we can introduce 
on the moduli space a K\"ahler metric, 
so-called the Weil-Peterson metric \cite{Tian},
through the K\"ahler potential $K(\psi,\bar\psi)=-\log i\int_M \Omega(\psi)
\wedge \bar\Omega(\psi)$ $= -\log i\sum_a (z_a(\psi)
\overline{ {\cal G}_a(\psi)}$ 
$- {\cal G}_a(\psi) \overline{ z_a(\psi)} )$. 
It is shown in ref.\cite{Str} that the prepotential 
$F(\psi)=$ 
${1\over 2}\sum_a z_a(\psi){\cal G}_a(\psi)$ describes the potential 
$K(\psi,\bar\psi)$ by 
\begin{equation}
K(\psi,\bar\psi)=
-{\rm log}i\sum_a \( z_a(\psi)\overline{{\pd F(\psi)\over \pd z_a}}
-\overline{z_a(\psi)} {\pd F(\psi) \over \pd z_a} \) \;,
\end{equation}
and defines the {\it special K\"ahler geometry} on the moduli space.  

Our definition (\ref{eqn: Fd}) of the prepotential, up to the prefactor 
$(w^{(0)})^{-2}$ which makes the prepotential invariant under 
$\Omega(\psi)\mapsto f(\psi)\Omega(\psi)$, implicitly contains a claim that 
$(w^{(0)},w^{(1)},{1\over2!}w^{(2)},-{1\over3!}w^{(3)}-
{c_2(X_{\D_s}) \over 12} w^{(1)})$ form the period integrals for a symplectic 
basis of $H_3(X_{\Delta_s^*},\ZZ)$. Several evidences for this claim are 
reported in \cite{HLY3}. 
In the following, we will restrict our attention to the form of 
the prepotential near a LCSL assuming its application to the enumerative 
geometry (the instanton counting).

\vspace{0.3cm}
Now consider the following expansion in the Chow ring associated to 
the series 
$w_0(x_\tau,0,A_\tau)=\sum_{n\in \ZZ_{\geq0}^{p-d}}c(n)x_\tau^n$;
\begin{equation}
\sum_{n\in \ZZ^{p-d}_{\geq0}}c\(n+{J\over 2\pi i}\)x_\tau^n=
\tilde w^{(0)}(x_\tau,J)+\tilde w^{(1)}(x_\tau,J)+
{1\over2!}\tilde w^{(2)}(x_\tau,J)
+{1\over3!}\tilde w^{(3)}(x_\tau,J) \;. 
\label{eqn: wJexpp}
\end{equation}

\vspace{0.3cm}\noindent
{\bf Lemma 6.3.} {\it 
The two definitions of the series (\ref{eqn: wJexp}) and (\ref{eqn: wJexpp}) 
are related by 
\begin{equation}
\begin{array}{crl}
w^{(0)}&=&\tilde w^{(0)} \;, \\
w^{(1)}&=&(\log x\cdt \hat J) w^{(0)} + \tilde w^{(1)} \;,\\
w^{(2)}&=&(\log x\cdt \hat J)^2 w^{(0)} + 
         2(\log x\cdt \hat J)\tilde w^{(1)} + \tilde w^{(2)}\;,\\
w^{(3)}&=&(\log x\cdt \hat J)^3 w^{(0)} + 
         3(\log x\cdt \hat J)^2\tilde w^{(1)}  +
         3(\log x\cdt \hat J) \tilde w^{(2)} + 
         \tilde w^{(3)}\;,\\
\end{array}
\end{equation}
where we have introduced an abbreviation $\log x\cdt\hat J=
\sum_{a=1}^{p-d} (\log x_a){1\over 2\pi i}J_a$.  }

\vspace{0.3cm}\noindent
{\bf Lemma 6.4.} {\it  
The series $\tilde w^{(d)}(x_\tau,J) \;\;(d=1,2,3)$ in (\ref{eqn: wJexpp}) 
have the form 
\begin{equation}
\begin{array}{crl}
\tilde w^{(1)}(x_\tau,J)&=&\sum_n c(n) \Psi^{(1)}(n) x_\tau^n \;\;, \\
\tilde w^{(2)}(x_\tau,J)&=&\sum_n c(n) 
          \{(\Psi^{(1)}(n))^2+\Psi^{(2)}(n)\}x_\tau^n \;\;, \\
\tilde w^{(3)}(x_\tau,J)&=&\sum_n c(n) \{(\Psi^{(1)}(n))^3+
     3\Psi^{(1)}(n)\Psi^{(2)}(n)+\Psi^{(3)}(n)\} x_\tau^n \;\;, \\
\end{array}
\label{eqn: wPsi}
\end{equation}
where $\Psi^{(k)}(n)$'s are elements in the Chow ring of degree $k$ 
defined by 
\begin{equation}
\begin{array}{crl}
&&\Psi^{(1)}(n)
 =-(\Jl_0)\psi(1-\nl_0)-\sum_{i=1}^{p} (\Jl_i) \psi(1+\nl_i) \;,\\
&&\Psi^{(2)}(n)=(\Jl_0)^2\psi'(1-\nl_0)-
     \sum_{i=1}^{p}(\Jl_i)^2\psi'(1+\nl_i) \;, \\
&&\Psi^{(3)}(n)=-(\Jl_0)^3\psi''(1-\nl_0)-
 \sum_{i=1}^{p}(\Jl_i)^3\psi''(1+\nl_i) \;, \\
\end{array}
\end{equation}
with $\hat J\cdt l_k = 
\sum_{a=1}^{p-d} {J_a \over 2\pi i} l^{(a)}_k \; (k=0,\cdots, p)$, 
$\psi(z)={d \; \over dz}\log\Gamma(z)$ and the derivatives of $\psi(z)$. }

\vspace{0.3cm}\noindent
{\bf Lemma 6.5.} {\it 
\begin{equation}
\Psi^{(1)}(0)=0 \;\;,\;\;
\Psi^{(2)}(0)=-{c_2(X_{\D_s})\over 12} \;\;,\;\; 
\Psi^{(3)}(0)=-{6\zeta(3) \over (2\pi i)^3} c_3(X_{\D_s}) \;\;. 
\end{equation}  }
\par
\noindent
{\bf (Proof)} 
These constant terms originate from those of the $\psi$-functions; 
$\psi(1)=-\gamma \;,\; \psi'(1)={\pi^2 \over 6} \;,\; 
\psi''(1)=-2 \zeta(3)$. These values of the $\psi$-functions combined 
with the adjunction formula for the total Chern class, with 
$D_i=J\cdt l_i$ under the rational equivalence in the Chow ring, 
\begin{equation}
c(X_{\D_s})={\prod_{i=1}^p(1+D_i) \over 1+[X_{\D_s}] }
           ={\prod_{i=1}^p(1+J\cdt l_i) \over 1-J\cdt l_0} \;, 
\end{equation}
result in our claim for the leading terms. 
(Note that $c_1(X_{\D_s})=0$ for $\Psi^{(1)}(0)$.)  \hfill $\Box$

\vspace{0.3cm}\noindent
{\bf Remark.} 
We can subtract these constant terms $\Psi^{(k)}(0)$ in a systematic way 
modifying the expansion (\ref{eqn: wJexpp}) slightly as follows;
\begin{equation}
\sum_{n\in \ZZ^{p-d}_{\geq 0}} 
{c(n+{J\over2\pi i}) \over c({J\over2\pi i})} x_\tau^n 
= w^{(0)}(x_\tau)+\tilde w_r^{(1)}(x_\tau,J)+
{1\over 2}\tilde w_r^{(2)}(x_\tau,J)
+ {1\over 3!} \tilde w_r^{(3)}(x_\tau,J) \;.
\end{equation}
This is because this change of normalization in the series $w_0$  
simply results in the replacement $\Psi^{(k)}(n)$ with 
$\Psi^{(k)}_r(n):= \Psi^{(k)}(n)-\Psi^{(k)}(0)$ in (\ref{eqn: wPsi}). 

\vspace{0.3cm}
Now it is immediate from  Lemmas 6.4 and 6.5 to obtain 

\vspace{0.3cm}\noindent
{\bf Lemma 6.6.} 
\begin{eqnarray*}
\tilde w^{(1)}&=&\tilde w^{(1)}_r \;\;,\;\;
\tilde w^{(2)}=-{c_2(X_{\D_s}) \over 12}w^{(0)} + \tilde w_r^{(2)} \;,\\
\tilde w^{(3)}&=&-{6\zeta(3) \over (2\pi i)^3} c_3(X_{\D_s})w^{(0)}
-{c_2(X_{\D_s}) \over 4} \tilde w_r^{(1)} + \tilde w_r^{(3)} . \\
\end{eqnarray*}
	\vspace{-2.5cm}
	\begin{equation} 
	\label{eqn: wr} \end{equation}
	\par\vspace{0.8cm}\noindent 

\vspace{0.3cm}
Now using the results in Lemmas 6.3-6.6, we may arrive at our final 
form of the prepotential, see also \cite{Sti}, modulo the kernel of 
the integration $\int_{X_{\D_s}}$ in Def.6.2;

\vspace{0.3cm}\noindent
{\bf Proposition 6.7.} {\it 
The invariant form of the prepotential ${\cal F}(x,J)$ 
may be expressed by  
\begin{eqnarray}
{\cal F}(t,J)
&=&{1\over6}(t\cdt J)^3-{c_2(X_{\D_s}) \over 24}(t\cdt J)+
{\zeta(3) \over 2(2\pi i)^3}c_3(X_{\D_s})   \nonumber \\
&& \quad\quad -{1\over 2}\log \left( \sum_{n\in \ZZ_{\geq0}^{p-d}} 
{ c(n+{J\over 2\pi i}) \over c({J\over 2\pi i}) } x_\tau^n \right) \;\;, 
\nonumber \\
\end{eqnarray}
with the mirror map $x_\tau=x_\tau(q)$.  }

\vspace{0.3cm}\noindent
{\bf Claim 6.8.} {\it 
Three times derivatives of the prepotential give the 
instanton corrected Yukawa couplings; 
\begin{eqnarray}
&&K_{t_a t_b t_c}(t)={\pd^3 \;\; \over \pd t_a \pd t_b \pd t_c} F(t) 
\nonumber \\
&&= \int_{X_{\D_s}} J_aJ_bJ_c + 
\sum_{ \Gamma\in H_2(X_{\D_s},\ZZ) \atop 
       \Gamma\not=0 } (\Gamma\cdt J_a)(\Gamma\cdt J_b)
(\Gamma\cdt J_c) N(\Gamma) 
{{\rm e}^{2\pi i \Gamma\cdt (t\cdt J)} \over 1- 
{\rm e}^{2\pi i \Gamma\cdt (t\cdt J)}} \;,  \nonumber \\
\end{eqnarray} 
where $\Gamma\cdt J:= \int_\Gamma J$ and $N(\Gamma)$ counts the number of 
the rational curves of class $\Gamma$ on the Calabi-Yau manifolds $X_{\D_s}$.}

\vspace{0.3cm}\noindent
{\bf Note.} 
Since the mirror map has the $q$-expansion 
$x_\tau^{(k)}(q)=q_k(1+{\cal O}(q))$, 
it is immediate to deduce that the number of lines $N(\Gamma)$ in 
$X_{\D_s}$ is counted by
\begin{equation}
N(\Gamma)=\int_{X_{\D_s}} -{1\over2}{c((\Gamma\cdt J)+J) \over c(J)} \;\;. 
\end{equation}
We see that the famous number $2785$ for the quintic in $\PP^4$ \cite{CdGP} 
is counted by this formula as 
\begin{equation}
N(1)=-{1\over2}\int_{\PP^4} 5J {(5+5J)(4+5J)(3+5J)(2+5J)(1+5J) 
\over (1+5J)^5}\;\;. 
\end{equation}

The invariant form of the prepotential ${\cal F}(x,J)$ may have significant
applications to extracting the predicted numbers of the rational curves
$N(\Gamma)$. In a recent work \cite{HSS}, this form has been utilized
efficiently to verify that the numbers $N(\Gamma)$ of a certain Calabi-Yau
manifold (Schoen's Calabi-Yau manifold) are related to the modular forms,
the theta function of the $E_8$ lattice and Dedekind's eta function. 


\def\thebibliography#1{\vskip 1.2pc{\centerline {\bf References}}\vskip 4pt
\list
 {[\arabic{enumi}]}{\settowidth\labelwidth{[#1]}\leftmargin\labelwidth
 \advance\leftmargin\labelsep
 \usecounter{enumi}}
 \def\newblock{\hskip .11em plus .33em minus .07em}
 \sloppy\clubpenalty4000\widowpenalty4000
 \sfcode`\.=1000\relax}
\let\endthebibliography=\endlist


\vskip 2pc
\noindent Shinobu Hosono\\
\noindent Department of Mathematics, Toyama University \\
\noindent Gofuku 3190, Toyama 930, Japan\\
\noindent email:hosono@sci.toyama-u.ac.jp 
\vskip 6pt

\end{document}